\documentclass[a4paper,12pt]{article}
\makeatletter
\@addtoreset{footnote}{page}
\makeatother

\usepackage{style-enumitem}

\usepackage{amsthm}
\usepackage{amsmath,amssymb,latexsym,amsfonts,mathrsfs}

\newcommand{\n}{\nonumber}

\renewcommand{\tilde}{\widetilde}
\renewcommand{\bar}{\overline}
\newcommand{\un}[1]{\underline{#1}}

\newcommand{\absol}[1]{\left| #1 \right|} %absolute value
\newcommand{\norm}[1]{\left\| #1 \right\|} %norm
\newcommand{\rbra}[1]{\!\left( #1 \right)} %round brackets or parentheses
\newcommand{\cbra}[1]{\!\left\{ #1 \right\}} %curly brackets or braces
\newcommand{\bD}{\ensuremath{\mathbb{D}}}
\newcommand{\bE}{\ensuremath{\mathbb{E}}}

\newcommand{\bN}{\ensuremath{\mathbb{N}}}

\newcommand{\bP}{\ensuremath{\mathbb{P}}}

\newcommand{\bR}{\ensuremath{\mathbb{R}}}

\theoremstyle{plain}
\newtheorem{Thm}{Theorem}[section]

\newtheorem{Lem}[Thm]{Lemma}

\newtheorem{Cor}[Thm]{Corollary}

\theoremstyle{definition}

\newtheorem{Def}[Thm]{Definition}
\newtheorem{Rem}[Thm]{Remark}

\newcommand{\Proof}[2][Proof]{\begin{proof}[{#1}] #2 \end{proof}}

\numberwithin{equation}{section}

\makeatletter
\renewcommand\section{\@startsection {section}{1}{\z@}%
                                   {-3.5ex \@plus -1ex \@minus -.2ex}%
                                   {2.3ex \@plus.2ex}%
                                   {\normalfont\large\bf}}
\makeatother

\makeatletter
\renewcommand\subsection{\@startsection {subsection}{1}{\z@}%
                                   {-3.5ex \@plus -1ex \@minus -.2ex}%
                                   {2.3ex \@plus.2ex}%
                                   {\normalfont\normalsize\bf}}
\makeatother

\renewcommand{\epsilon}{\varepsilon}

\usepackage{color}
\usepackage{mathtools}
\mathtoolsset{showonlyrefs=true}
\begin{document}

\begin{center}
{\Large \bf 
Generalized refracted L\'{e}vy process and its application to exit problem
}
\end{center}
\begin{center}
Kei Noba\footnote{Department of Mathematics, Graduate School of Science, Kyoto University
Sakyo-ku, Kyoto 606-8502, Japan. Email: knoba@math.kyoto-u.ac.jp (K. Noba), kyano@math.kyoto-u.ac.jp (K. Yano).}
and Kouji Yano$^\ast$
\end{center}

\begin{abstract}
Generalizing Kyprianou--Loeffen's refracted L\'evy processes, we define a new refracted 
L\'evy process which is a Markov process whose positive 
and negative motions are L\'evy processes 
different from each other. 
To construct it we utilize the excursion theory. 
We study its exit problem and the potential measures of the 
killed processes. 
We also discuss approximation problem. 
\end{abstract}
\

%%%%% text %%%%%

\section{Introduction}
\emph{Exit problem} of a real-valued stochastic process $Z=\{Z_t : t \geq 0\}$ 
is the problem to characterize the law of the first time of exiting an interval $[b, a]$ 
for $b<a $. 
In this paper, we are interested in the Laplace transform 
\begin{align}
\bE^Z_x\rbra{ e^{-q\tau_{a}^{+}} ;\tau_{a}^{+} <\tau_{b}^{-}}
\label{101h}
\end{align}
for $q\geq0$ and a starting point $x\in[b, a]$, where 
\begin{align}
\tau_{a}^{+} = \inf\{ t>0: Z_t > a\} \text{ and } 
\tau_{b}^{-} = \inf\{t>0 : Z_t < b \}.
\label{hitting time}
\end{align}
When $Z$ is a spectrally negative L\'{e}vy process, it is well known that 
\begin{align}
\bE^Z_x\rbra{ e^{-q\tau_{a}^{+}} ;\tau_{a}^{+} <\tau_{b}^{-}}
=\frac{W_Z^{(q)}(x-b)}{W_Z^{(q)}(a-b)}, 
\end{align}
where $W_Z^{(q)}$ is the $q$-scale function of $Z$. 
\par
Kyprianou and Loeffen \cite{Kyp2} have studied the exit problem when $Z$ was 
a \emph{refracted L\'{e}vy process} $U$, which 
was defined as the strong solution of 
the stochastic differential equation
\begin{align}
U_t- U_0=X_t-X_0+\alpha\int_0^t1_{\{U_s < 0\}}ds~~~t \geq0,
\label{Kyprianou--Loeffen's refracted}
\end{align}
where the driving noise $X$ is a spectrally negative L\'evy process 
and $\alpha$ is a positive constant. 
Define $Y_t=X_t+\alpha t$. 
Then the positive and negative motions of $U$ is given as 
\begin{align}
U_t -U_s= 
\begin{cases}
X_t - X_s \text{ whenever } U_r \geq 0 \text{ for any } r \in [s, t) \\
Y_t - Y_s ~\text{ whenever } U_r < 0 \text{ for any } r \in [s, t). 
\end{cases}
\end{align} 
They proved that 
the Laplace transform \eqref{101h} for $Z=U$ takes the form
\begin{align}
\bE^U_x\rbra{ e^{-q\tau_{a}^{+}} ;\tau_{a}^{+} <\tau_{b}^{-}}
=\frac{W_U^{(q)} (x,  b)}{W_U^{(q)} (a, b)}, 
\label{1.05q}
\end{align}
where the function $W_U^{(q)}$ is defined by 
\begin{align}
W_U^{(q)} (x, y) = W_{Y}^{( q )} \rbra{x-y}
	+ \alpha 1_{(x \geq 0)} \int ^{x}_{0} W_{X}^{(q)}
	(x-z) {W_{Y}^{(q) \prime}} (z-y)dz \label{1.06q}
\end{align}
with $W_Y^{(q)}$ being the $q$-scale function of $Y$. 
They obtained, in addition, a representation of the potential measures of $U$ 
using scale functions of $X$ and $Y$. 
\par
A spectrally negative L\'{e}vy process can be regarded as the capital 
of an insurance company and applied to evaluate the risk of ruin. Hence it is 
sometimes called a \emph{L\'{e}vy insurance risk process}.
The Kyprianou--Loeffen's refracted L\'evy process $U$ can be regarded as 
a modified insurance risk process when dividends are being  paid out at a rate $\alpha$ 
during the period it exceeds $0$. 
 \par
In this paper, we generalize Kyprianou--Loeffen's refracted L\'{e}vy processes. 
For two L\'{e}vy processes $X$ and $Y$ which may have 
different L\'{e}vy {exponents}, we construct a new refracted process 
whose positive and negative motions have the same law as $X$ and $Y$, respectively. 
More precisely, 
\begin{align}
\begin{cases}
\text{If }x>0 \text{, } \rbra{U_t}_{t \leq \tau^-_0} 
			\text{ under }\bP^U_x \text{ is equal in law to } 
		\rbra{X_t}_{t \leq \tau^-_0} \text{ under }\bP^X_x\\
\text{If }x<0 \text{, } \rbra{U_t}_{t \leq \tau^+_0} 
			\text{ under }\bP^U_x \text{ is equal in law to } 
		\rbra{Y_t}_{t \leq \tau^+_0} \text{ under }\bP^Y_x.
\end{cases}
\end{align}
One may expect that we can characterize the desired process as a solution to the 
following stochastic differential equation
	\begin{align}
	\label{ref general}
	U_t -U_0=  \int_{(0, t]} 1_{\{U_{s-}  \geq 0\}} d X_{s}
	+ \int_{(0, t]} 1_{\{U_{s-} < 0\} } d Y_{s}, 
	\end{align}
where the driving noises $X$ and $Y$ are supposed to be independent. 
Although \eqref{ref general} for $Y_t\overset{d}{=} X_t + \alpha t$ 
is apparently different from \eqref{Kyprianou--Loeffen's refracted} because of independence, 
their solutions are actually equivalent in law. 
When $X$ has bounded variation paths, we can construct a solution of \eqref{ref general} 
by a simple method of piecing excursions (see \cite{Kyp2}); 
otherwise we do not know existence of a solution of \eqref{ref general}. 
When $X$ and $Y$ are compound 
Poisson processes with positive drifts, uniqueness of the solution is easily proved 
because of the fact that the point $0$ is irregular for itself for any solution $U$; 
otherwise we do not know uniqueness of a solution of \eqref{ref general}. 
\par
In this paper we utilize the excursion theory instead of a stochastic differential 
equation. 
Let $X$ and $Y$ be two spectrally negative L\'evy processes. 
Suppose $X$ has unbounded variation paths and 
has no Gaussian component. 
We then define the excursion measure $n^U$ by 
	\begin{align}
&n^{U}\rbra{ F\rbra{(U_t)_{t < \tau^{-}_{0}}, (U_{t+\tau^{-}_{0}})_{t\geq 0}}}
=n^{X}\rbra{ \mathbb{E}^{Y^0}_{y}
\rbra{F\rbra{w, (Y^0_t)_{t\geq 0}}}
{\biggr|}_{\tiny{\begin{subarray} xy=X(\tau^{-}_{0})\\w=(X(t))_{t<\tau^{-}_{0}}\end{subarray}}}},
\label{excursion measure of U}
	\end{align}
where $n^X$ stands for an excursion measure of $X$ 
and $Y^0_t=Y_{t \land T_0}$ for the stopped process of $Y$ upon hitting zero. 
We define the stopped process $\bP^{U^0}_x$ by 
\eqref{excursion measure of U} with $n^X$ being replaced by $\bP^{X}_x$. 
We can therefore construct a Feller process from $n^U$ together with 
the family of stopped processes $\cbra{\bP^{U^0}_x}_{x \neq 0}$. 
As one of our main theorems, we show the Laplace transform \eqref{101h} for 
the process $Z=U$, 
our new refracted L\'evy process, takes the same form as \eqref{1.05q} 
where $W_U^{(q)}$ will be defined in 
{Theorem \ref{refracted exit problem}} 
in a more complicated form than \eqref{1.06q}.
Note that $W_U^{(q)}$'s will be represented using only Laplace exponents and 
scale functions of $X$ and $Y$. 
Furthermore, we will study the potential measures of $U$ 
with and without absorbing barriers. 
\par
We finally discuss approximation problem. 
Let $X$ and $Y$ be as in the previous paragraph. 
Let $X^{(n)}$ and $Y^{(n)}$ be the compound Poisson processes with positive drifts 
obtained from $X$ and $Y$, respectively, 
by removing small jumps of magnitude less than $\frac{1}{n}$. 
Assuming that $X^{(n)}$ and $Y^{(n)}$ are independent, we construct $U^{(n)}$ 
as the unique solution of \eqref{ref general}. We thus show that $U^{(n)}$ 
converges to our refracted process $U$ in law on the space of c\`adl\`ag paths 
equipped with the Skorokhod topology. 
\par
The organization of the present paper is as follows. 
In Section 2 we propose some notation and recall preliminary facts about 
spectrally negative L\'evy processes.  
In Section 3 we calculate several quantities
related to excursion measures and scale functions.
In Section 4 we recall Kyprianou--Loeffen's refracted L\'evy processes.
In Section 5 we define our new refracted L\'evy processes. 
In Section 6 we study the exit problem of our refracted L\'evy processes. 
In Section 7 we calculate the potential measures of our 
refracted L\'evy processes killed upon exiting $[b, a]$.
In Section 8 we study the approximation problem. 
{In Section \ref{construction} 
we make a careful treatment of Markov property of our new process. }

\section{Notation and preliminaries}
Let $\bD$ denote the set of functions 
	$ \omega : [0, \infty ) \rightarrow \bR$ which are c\`adl\`ag. 
We equip $\bD$ with the Skorokhod topology.
Let ${\cal{B}}(\bD)$ denote the class of Borel sets of $\bD$. 
\par
When we consider a process $Z=\{Z_t : t \geq 0\}=\{Z(t) : t \geq 0\}$, 
we always write $\bP^Z_x$ for 
the underlying probability measure for $Z$ starting from $x$. 
In addition to the passage times $\tau^+_a$ and $\tau^-_b$ defined in \eqref{hitting time}, 
we sometimes need the hitting time of a point $x\in \bR$ denoted by 
\begin{align}
T_x=\inf\{ t>0:Z_t=x\}.
\end{align}
For $q>0$, $x \in \bR$ and a non-negative or bounded measurable function $f$, 
we write 
		\begin{align}\label{}
				R^{(q)}_{Z}f(x) := \mathbb{E}^{Z}_x \rbra{\int_{0}^{\infty} e^{-qt} 
		 f(Z_t)  dt }.
		\end{align}
	We write $r_Z^{(q)}(x, y)$ for the resolvent density, if it exists, i.e., 
		\begin{align}
		R^{(q)}_{Z}f(x)=\int_{\bR} r^{(q)}_{Z}\rbra{x, y} f(y) dy .
		\label{202p}
		\end{align}
		We sometimes settle a lower barrier $b<0$ and an upper barrier $a>0$. 
		For $q>0$, $ x \in \bR$ and a non-negative or bounded 
		measurable function $f$, we write 
		\begin{align}\label{}
		\underline{\overline{R}}^{(q;b,a)}_{Z}f(x):=
				\underline{\overline{R}}^{(q)}_{Z}f(x) := \mathbb{E}^{Z}_x 
		\rbra{\int_{0}^{\tau^{+}_{a}\land \tau^{-}_{b}} e^{-qt}  f(Z_t)  dt }, 
		\\
		\underline{R}^{(q; b)}_{Z}f(x):=
		\underline{R}^{(q)}_{Z}f(x) := \mathbb{E}^{Z}_x 
		\rbra{\int_{0}^{\tau^{-}_{b}} e^{-qt}  f(Z_t)  dt },
		\\
		\overline{R}^{(q; a)}_{Z}f(x):=
		\overline{R}^{(q)}_{Z}f(x) := \mathbb{E}^{Z}_x 
		\rbra{\int_{0}^{\tau^{+}_{a}} e^{-qt}  f(Z_t)  dt }.
		\end{align}
	We write 
	$\underline{\overline{r}}^{(q; b, a)}_{Z}\rbra{ x, y}=
	\underline{\overline{r}}^{(q)}_{Z}\rbra{ x, y}$, 
	{$\underline{{r}}^{(q; b, a)}_{Z}\rbra{ x, y}=
	\underline{{r}}^{(q)}_{Z}\rbra{ x, y}$ 
	and 
	${\overline{r}}^{(q; b, a)}_{Z}\rbra{ x, y}=
	{\overline{r}}^{(q)}_{Z}\rbra{ x, y}$ 
	for the corresponding densities, if they exist.} 
\par
Let $Z$ 
%$\rbra{ Z=\{ Z_t : t \geq 0 \} , \mathbb{P}^{Z}_{x}}$ 
be a spectrally negative L\'{e}vy process, 
which is always assumed not to be monotone. 
Then it is well known that the Laplace exponent
	\begin{align}\label{}
	\Psi_{Z} (q) := \log \mathbb{E}^{Z}_{0} \rbra{e^{q Z_{1}}} 
	\end{align}
is finite for all $q \geq 0$.
We denote its right inverse by
	\begin{align}\label{}
	\Phi_{Z} (\theta) 
	= \inf \{ q  \geq 0 : \Psi_{Z} \rbra{q} = \theta \},
	\end{align}
which is finite for all $\theta \geq 0$.
If $Z$ has bounded variation paths,
the Laplace exponent is known to necessarily take the form
	\begin{align}\label{208x}
	\Psi_{Z}(q)
	=\delta_{Z} q - \int_{(- \infty , 0)}
	\rbra{1-e^{q y}} \Pi_{Z} (dy)
	\end{align} 
for some constant $ \delta_Z > 0$ and some L\'{e}vy measure $\Pi_{Z}$ satisfying 
$\Pi_{Z}[ 0, \infty ) = 0$ and 
$\int_{(-\infty , 0)}\rbra{1\land |y|}\Pi_Z (dy) < \infty$. 
If $Z$ has unbounded variation paths,
the Laplace exponent is known to necessarily take the form 
	\begin{align}
	\Psi_Z(q) 
	= \gamma_Z q + \frac{\sigma_{Z}^{2}}{2} q^{2} - 
	\int_{(-\infty , 0)}\rbra{ 1 - e^{q y} + 
	q y 1_{(-1, 0)}(y) } \Pi_Z (dy)
	\label{Laplace exponent}
	\end{align}
for some constants $\gamma_Z \in \mathbb{R}$ and $\sigma_Z \geq 0$ and 
some L\'{e}vy measure $\Pi_{Z}$ 
satisfying $\Pi_{Z}[0, \infty ) = 0$ and 
$\int_{(-\infty , 0)} \rbra{1 \land y^{2}} \Pi_Z (dy) < \infty$.
	\begin{Def}
	For each $q \geq 0$, we define $W_{Z}^{(q)}: \bR \rightarrow [0, \infty)$ such that 
	 $W_{Z}^{(q)}= 0$ on $(-\infty, 0)$ and $W_{Z}^{(q)}$ on $[0, \infty )$ is continuous 
	satisfying
		\begin{align}
		\int_{0}^{ \infty} e^{-\beta x}W_{Z}^{(q)}(x)dx =\frac{1}{\Psi_{Z}(\beta) - q}
		\end{align}
	for all $\beta > \Phi_{Z}(q)$. This function $W_Z^{(q)}$ is called the 
	{\emph{$q$-scale function of $Z$}}.
	\end{Def}
	For the proof of unique existence and its basic facts listed below, see, e.g., \cite{Kyp}. 
	For all $b < x< a$ and $q \geq 0$, we have
		\begin{align}
		\bE^{Z}_{x}\rbra{e^{-q\tau^{+}_{a}};\tau^{+}_{a}<\tau^{-}_{b} }
		= \frac{W_{Z}^{(q)}(x-b)}{W_{Z}^{(q)}(a-b)}
		\label{exit problem}
		\end{align}
	and
	\begin{align}
	\bE^{Z}_{x}\rbra{e^{-q\tau^{+}_{a}};\tau^{+}_{a}<\infty }
		=e^{-\Phi_{Z}(q)(a-x)}.
	\label{one-sided exit problem}
	\end{align}
	It is known that, when $Z$ has bounded variation paths, 
	we have 
	\begin{align}\label{scale function 0}
	W_Z^{(q)}(0) = \frac{1}{\delta_Z}
	\end{align}
	for all $q \geq 0$.
	For all $q \geq 0$, we have
	\begin{align}
	r^{(q)}_{Z}(x, y)=\Phi_{Z}'(q)e^{-\Phi_{Z}(q)(y-x)}-W_{Z}^{(q)}(x-y),~~~~x, y\in \bR,
	\label{potential density}
	\end{align}
	\begin{align}
	\underline{r}_{Z}^{(q; b)}( x, y) =\underline{r}_{Z}^{(q)}( x, y) 
		= e^{-\Phi_{Z}(q)(y-b)}W_{Z}^{(q)}(x-b) - W_{Z}^{(q)}(x-y),~~ x , y \in [b, \infty),
	\label{one-sided potential density}
	\end{align}
	\begin{align}
	\overline{r}_{Z}^{(q; a)}( x, y)=\overline{r}_{Z}^{(q)}( x, y) 
	= e^{-\Phi_{Z}(q)(a-x)}W_{Z}^{(q)}(a-y) - W_{Z}^{(q)}(x-y),
	~~x , y\in (-\infty, a]	
	\label{one-sided potential density 2}
	\end{align}
	and
	\begin{align}
	\underline{\overline{r}}_{Z}^{(q; b,a)}( x, y)=
	\underline{\overline{r}}_{Z}^{(q)}( x, y) = 
		\frac{W_{Z}^{(q)}(x - b)W_{Z}^{(q)}(a-y)}{W_{Z}^{(q)}(a-b)}-W_{Z}^{(q)}(x-y),~~
		x , y\in [b, a].
		\label{killed potential density}
	\end{align}
	We write {$\tilde{\Pi}_Z$} for the measure carried on $(-\infty, 0)\times(0, \infty)$ 
	defined by 
	\begin{align}
	\tilde{\Pi}_{Z}(du~ dv) : = {1_{\{u<0 , \, v>0\}}}\Pi_{Z}(du -v)dv.
	\end{align}
	\begin{Thm}\label{Gerber-Shiu}
	{$(i)$} 
	For all $0<x<\infty$, $q \geq 0$, and non-negative measurable function 
	$f:\bR^{2} \rightarrow [0, \infty)$, we have
		\begin{align}
		&\bE^{Z}_{x}\rbra{e^{-q\tau^{-}_{0}}f(Z_{\tau^{-}_{0}}, Z_{\tau^{-}_{0}-})
			; \tau^{-}_{0} < \infty, Z_{\tau^{-}_{0}} < 0}
		=\int f(u, v) G^{(q)}_{Z}(x, du ~dv),
		\end{align}
	where $G^{(q)}_{Z}(x, \cdot)$ is the measure carried on 
	$(-\infty, 0) \times (0, \infty)$ defined by 
		\begin{align}
		G^{(q)}_{Z}(x, du~ dv) 
		:= {\underline{r}}_Z^{(q; 0)}(x, v) \tilde{\Pi}_{Z}(du~ dv).
		\end{align}
	\par
	{$(ii)$} 
	For all $0<x<a$, $q \geq 0$, and non-negative measurable function $f$, we have
		\begin{align}
		&\bE^{Z}_{x}\rbra{e^{-q\tau^{-}_{0}}f(Z_{\tau^{-}_{0}}, Z_{\tau^{-}_{0}-})
			; \tau^{-}_{0} < \tau^{+}_{a}, Z_{\tau^{-}_{0}} < 0}
		=\int f(u, v)\overline{G}_{Z}^{(q, a)}(x, du ~dv),
		\label{219q}
		\end{align}
	where $\overline{G}^{(q, a)}_{Z}(x, \cdot)$ is the measure carried on 
	$(-\infty, 0) \times (0, \infty)$ defined by 
		\begin{align}
		\overline{G}_{Z}^{(q, a)}(x,  du~ dv)=\overline{G}_{Z}^{(q)}(x,  du~ dv)
		:= \overline{\underline{r}}_Z^{(q;0,a)}(x, v)\tilde{\Pi}_{Z}(du~ dv).
		\label{220q}
		\end{align}
	\end{Thm}
	{We omit the proof of Theorem \ref{Gerber-Shiu} because it can be found in 
	\cite[Theorem 10.1 and Exercise 10.6]{Kyp} and also in 
	\cite[Theorem 5.5]{Kyp2}. 
	These kernels $G_Z^{(q)}$ and $\overline{G}_{Z}^{(q)}$ 
	are called the \emph{Gerber--Shiu measures}. }

\section{Some calculations related to excursion measures and scale functions}

In this section, we make some calculations related to excursion measures 
and scale functions for a spectrally negative L\'evy process $X$. 
	See \cite{Par}, \cite{Avr} and \cite{Par2} for recent studies on 
	a close relation between $n^X$, {i.e.,} the excursion measure of $X$ itself, 
	and the excursion measure of the reflected of $X$. 
\par
We divide the discussion into the two cases of unbounded and of bounded variations. 
\par
{\bf{(I)}} We assume that 
\begin{align}
\text{$X$ has unbounded variation paths and have no Gaussian component}. 
\end{align}
Since $0$ is regular for itself, $X$ has an excursion measure $n^X$ away from zero. 
We impose on $n^{X}$ the following normalization:
	\begin{align}
	n^{X}\rbra{ 1- e^{-qT_0}}
	=\frac{1}{r^{(q)}_{X}(0, 0)}
	=\frac{1}{\Phi_{X}'(q)}
	=\Psi_{X}'\rbra{ \Phi_{X}(q)} .
	\label{normalization of excursion}
	\end{align} 
Note that $n^X$ is carried on the set of c\`adl\`ag paths stopped upon hitting $0$. 
Note also that $n^X$ possesses the Markov property; for example, 
	\begin{align}
	n^X\rbra{ X_s \in B_1, X_t \in B_2}
	=n^X\rbra{ 1_{\{ X_s \in B_1\}} \bP^{X^0}_{X(s)}\rbra{X^0_{t-s} \in B_2}}, 
	\end{align} 
for all $0< s < t$ and $B_1$, $B_2 \in {\cal{B}}(\bR)$, 
where $X^0_t=X_{t\land T_0}$ denotes the stopped process of $X$ upon hitting zero. 
Since $X$ has no Gaussian component, we can see 
\begin{align}
0<\tau^{-}_{0}<T_0 \leq \infty ~ \text{ or } ~  \tau^{-}_{0}= T_0 =\infty ~~~~~n^X\text{-a.e.}
\end{align} 
by \cite[Theorem 3]{Par}. 
%For all $t\geq 0$, we denote 	$\overline{X}=\sup_{s\geq0}X_s$ and 
%	$\underline{X}=\inf_{s \geq 0}X_s$. 
	\begin{Thm}\label{equation lemma}
	For all $a > 0$ and $q \geq 0$, we have 
		\begin{align}
		n^{X}\rbra{e^{-q\tau^{+}_{a}};\tau^+_a < \infty}=\frac{1}{W_{X}^{(q)}(a)}.
		\label{303p}
		\end{align}
	In particular, we have
		\begin{align}
		n^{X}\rbra{\tau^+_a < \infty}=\frac{1}{W_{X}(a)}, 
		\label{306pika}
		\end{align}
	where $W_{X} := W^{(0)}_{X}$.
	\end{Thm}
	\begin{Rem}
	The two ways of normalization \eqref{normalization of excursion} and \eqref{306pika} 
	are natural analogies of those for diffusion processes. 
	See %e.g. 
		\cite[$(2.5)$ and Theorem $3.1$]{Che} 
		and \cite[$(39)$ and Theorem $3.1$]{Yan}.
	\end{Rem}
	{\begin{Rem}\label{Rem303}
	The left-hand side of \eqref{303p} may admit several other {expressions}, such as 
	\begin{align}
	n^X\rbra{e^{-q\tau^+_a}}
	=n^X\rbra{e^{-q\tau^+_a} ; \tau^+_a < \infty}
	=n^X\rbra{e^{-q\tau^+_a} ; \tau^+_a < \tau^-_0 }. \label{307ac}
	\end{align}
	The first equality of \eqref{307ac} follows from the fact that $e^{-q\tau^+_a}=0$ 
	on $\{\tau^+_a=\infty\}$. 
	Since $X$ has no positive jumps, the measure $n^X$ is supported {on} 
	the disjoint union 
	\begin{align}
	\{\tau^+_a < \tau^-_0 \leq \infty \}\cup\{\tau^-_0<\tau^+_a = \infty\}\cup 
	\{\tau^-_0 = \tau^+_a = \infty \}. 
	\end{align}
	Thus the sets $\{\tau^+_a < \infty\}$ and $\{\tau^+_a< \tau^-_0\}$ are equal up to 
	$n^X$-null sets, which yields the second equality of \eqref{307ac}. 
	\end{Rem}
	}
	The following theorem can be regarded as the Gerber--Shiu measure for 
	the excursion measure (see also \cite{Par2}).
		\begin{Thm}
	\label{excursion Gerber-Shiu}
	For all $q\geq 0$ and non-negative measurable function $f$, we have
		\begin{align}
		n^{X}\rbra{
		e^{-q\tau^{-}_{0}}f(X_{\tau^{-}_{0}}, X_{\tau^{-}_{0}-}) ; \tau^{-}_{0} < \infty}
		=\int f(u, v) K^{(q)}_{X}(du~ dv), 
		\label{317p}
		\end{align} 
	where $K_X^{(q)}$ is the measure carried on $(-\infty, 0)\times (0, \infty)$ defined by 
		\begin{align}
		K^{(q)}_{X}(du~d v) 
		=  e^{-\Phi_{X}(q)v} \tilde{\Pi}_{X}(du~ dv). 
		\label{314aa}
		\end{align}
%	(The constant $c$ will turn out to be $1$ in Lemma \ref{constant}.)
	\end{Thm}
	We prove Theorems \ref{equation lemma} and \ref{excursion Gerber-Shiu} at the 
	same time. 
	\Proof[Proof of Theorems \ref{equation lemma} and \ref{excursion Gerber-Shiu}]{
	\ \ \par
	{\bf{Step.1}} We show that the quantity
		\begin{align}
		c:=n^{X}\rbra{e^{-q\tau^{+}_{a}};\tau^+_a < \infty}W_{X}^{(q)}(a)
		\end{align}
	does not depend upon $a > 0$ nor $q \geq 0$.
	First, we prove 
		\begin{align}
		\label{scale1}
		n^{X}\rbra{\tau^+_a < \infty}W_{X}(a)
		=&n^{X}\rbra{e^{-q\tau^{+}_{a}};\tau^+_a < \infty}
		W_{X}^{(q)}(a) 
		\end{align}
	for all $a > 0$ and $q \geq 0$. Using the monotone convergence theorem, we have
		\begin{align}
		\frac{n^{X}\rbra{e^{-q\tau^{+}_{a}};\tau^+_a < \infty}}
			{n^{X}\rbra{\tau^+_a < \infty}}
%		=&\lim_{\epsilon \downarrow 0}
%			\frac{n^{X}\rbra{e^{-q\tau^{+}_{a}};\tau^+_a < \infty}}
%			{n^{X}\rbra{e^{-q\tau^{+}_{\epsilon}};\tau^+_a < \infty}}\\
		=&\lim_{\epsilon \downarrow 0}
			\frac{n^{X}\rbra{e^{-q\tau^{+}_{\epsilon}}\rbra{e^{-q\tau^{+}_{a}} 
			1_{\{\tau^+_a < \infty\}}}
			\circ \theta_{\tau^{+}_{\epsilon}} ; \tau^+_\epsilon < \infty}}
			{n^{X}\rbra{e^{-q\tau^{+}_{\epsilon}}
			\rbra{ 1_{\{\tau^+_a < \infty\}}}\circ \theta_{\tau^{+}_{\epsilon}}
			; \tau^+_\epsilon < \infty }}.
			\label{305}
		\end{align}
	Using the strong Markov property and \eqref{exit problem}, we have
	\begin{align}
		\eqref{305}
		=&\lim_{\epsilon \downarrow 0}\frac{n^{X}\rbra{e^{-q\tau^{+}_{\epsilon}};
		\tau^+_\epsilon < \infty}
		\mathbb{E}^{X}_{\epsilon}\rbra{e^{-q\tau^{+}_{a}}1_{\{\tau^{+}_{a} 
		< \tau^{-}_{0} \}}}}
		{n^{X}\rbra{e^{-q\tau^{+}_{\epsilon}};\tau^+_\epsilon < \infty}
		\mathbb{P}^{X}_{\epsilon}\rbra{\tau^{+}_{a} 
		< \tau^{-}_{0} }}
%		=&\lim_{\epsilon \downarrow 0}\frac{
%		W_{X}(a)W_{X}^{(q)}(\epsilon)}
%		{W_{X}^{(q)}(a)W_{X}(\epsilon)}\\
		=\frac{W_{X}(a)}{W_{X}^{(q)}(a)},
		\label{308}
		\end{align}
	where we used 
	$\lim_{\epsilon \downarrow 0}W_X^{(q)}(\epsilon)/W_X(\epsilon) =1$ by 
	\cite[Lemma 1. (i)]{Don}.
	Second, we prove 
		\begin{align}
		n^{X}\rbra{\tau^+_{a_1} < \infty}W_{X}(a_1)
		=&n^{X}\rbra{\tau^+_{a_2} < \infty}W_{X}(a_2),
		\end{align}
	for all $0<a_1 < a_2$. This identity can be obtained by
		\begin{align}
		\frac{n^{X}(\tau^+_{a_2} < \infty)}{n^{X}(\tau^+_{a_1} < \infty)}
		=&\frac{n^{X}\rbra{
			\rbra{\tau_{a_2}^{+} <\tau_{0}^{-}} \circ \theta_{\tau^{+}_{a_1}}
			;\tau^+_{a_1} < \infty}}{n^{X}\rbra{\tau^+_{a_1} < \infty}}
%		\\
%		=&\frac{n^{X}\rbra{\tau^+_{a_1} < \infty}\bE^{X}_{a_1}
%			\rbra{\tau_{a_2}^{+} <\tau_{0}^{-}}}{n^{X}\rbra{\tau^+_{a_1} < \infty}}
%		\label{311}
%		\\
		=\frac{W_{X}(a_1)}{W_{X}(a_2)},
		\label{312}
		\end{align}
	where we used the strong Markov property and 
	\eqref{exit problem}. \par
	{\bf{Step.2}} We show \eqref{317p} with $K_X^{(q)}$ being multiplied by $c$. 
	Using the monotone convergence theorem and the strong Markov property, we have 
		\begin{align}
		&n^{X}\rbra{e^{-q\tau^{-}_{0}}f\rbra{X_{\tau^{-}_{0}}, X_{\tau^{-}_{0}-}}
			; \tau^{-}_{0} < \infty}\nonumber\\
%		=&\lim_{\epsilon \downarrow 0}n^{X}\rbra{
%		e^{-q\tau^{-}_{0}}f\rbra{X_{\tau^{-}_{0}}, X_{\tau^{-}_{0}-}}
%			; X_{\tau^{-}_{0}-} > \epsilon,  \tau^{-}_{0} < \infty}
%		\\
		=&\lim_{\epsilon \downarrow 0}n^{X}\rbra{e^{-q\tau^{+}_{\epsilon}}
			\rbra{e^{-q\tau^{-}_{0}}f\rbra{X_{\tau^{-}_{0}}, X_{\tau^{-}_{0}-}}
			1_{\cbra{ X_{\tau^{-}_{0}-} > \epsilon}}} \circ
			\theta_{\tau^{+}_{\epsilon}}; \tau^+_\epsilon < \infty,  \tau^{-}_{0} < \infty}
%		\\
%		=&\lim_{\epsilon \downarrow 0}
%			n^{X}\rbra{ e^{-q\tau^{+}_{\epsilon}};\tau^+_\epsilon < \infty }
%			\mathbb{E}^{X}_{\epsilon}
%			\rbra{ e^{-q\tau^{-}_{0}}f\rbra{X_{\tau^{-}_{0}}, X_{\tau^{-}_{0}-}}
%			; X_{\tau^{-}_{0}-} > \epsilon ,  \tau^{-}_{0} < \infty},
		\label{3018}
		\end{align}
		and using {$(i)$ of Theorem \ref{Gerber-Shiu}} 
		and {\bf{Step.1}}, we have
		\begin{align}
		\eqref{3018}%\nonumber
		=&\lim_{\epsilon \downarrow 0}
			\frac{c}{W^{(q)}_{X}(\epsilon)}
			W^{(q)}_{X}(\epsilon)
			\int_{\epsilon}^{\infty}dv
			\int_{(-\infty,0)}
			e^{-\Phi_{X}(q)v}f(u, v)
			\Pi_{X}(du - v)
		\\
		=& \int f(u, v) 
			c K^{(q)}_{X}(du, dv).
		\end{align}
	\par
%	which coincides with the right hand side of \eqref{317p}. 
	{\bf{Step.3}} We show $c =1$. 
	Since $X$ has no Gaussian component, i.e., $\sigma_X=0$, 
	differentiating \eqref{Laplace exponent}, we have
		\begin{align}
			\Psi_X'(q) 
		= \gamma_X  +
		\int_{(-\infty , 0)}\rbra{ ye^{q y} - y 1_{(-1, 0)} (y)} \Pi_X (dy)
		\label{derivative of Laplace exponent}
		\end{align}
	for all $q > 0$.
	Using \eqref{normalization of excursion}, we have on one hand 
		\begin{align}
		n^X\rbra{1-e^{-qT_0}}=&\Psi_{X}'\rbra{\Phi_{X}(q)}
		= \gamma_X  +
		\int_{(-\infty , 0)}
		\rbra{ ye^{\Phi_{X}(q) y} - y 1_{(-1, 0)}(y) } \Pi_X (dy). 
		\label{326b}
		\end{align}
	On the other hand, 
	using the monotone convergence theorem and the strong Markov property, we have
	\begin{align}
	&n^{X}\rbra{1 - e^{-qT_0}}\nonumber\\
		=&n^X\rbra{\tau^{-}_{0} = \infty}
		+n^{X}\rbra{1 - e^{-qT_0} ; \tau^{-}_0 < \infty}\\
		%=&n^X\rbra{\tau^+_1 < \infty}\bE^{X}_1\rbra{\tau^{-}_{0}= \infty}
		%+n^X\rbra{1 - e^{-q\tau_0^-}
		%	\rbra{e^{-qT_0}}\circ\theta_{\tau^-_0};\tau^{-}_{0} < \infty}\\
		=&n^X\rbra{\tau^+_1 < \infty}\lim_{p \uparrow \infty}
			\bE^{X}_{1}\rbra{\tau^{+}_{p}<\tau^-_0}
		+n^X
		 \rbra{1- e^{-q\tau^-_0}\bE^X_{X(\tau^-_0)}\rbra{e^{-qT_0}}; \tau^-_0 < \infty }.
		\label{30026}
	\end{align}
	Using \eqref{exit problem}, \eqref{one-sided exit problem}, 
	{\bf{Step.1}}  
	and {\bf{Step.2}}, we have
		\begin{align}
		\eqref{30026}
		=&n^{X}\rbra{\tau^+_1 < \infty}\lim_{p \uparrow \infty}\frac{W_{X}(1)}{W_{X}(p)}
		+n^X\rbra{1-e^{-q\tau^-_0+\Phi_{X}(q)X(\tau^-_0)};\tau^-_0 <\infty}\\
		=&c\frac{1}{W_X(\infty)}
		+c\int 
		\rbra{e^{-\Phi_{X}(0)v} - e^{\Phi_{X}(q)(u-v)}}{\tilde{\Pi}}_X(du ~dv).
		\label{328a}
		\end{align}
	Since it is known that
		\begin{align}
		W_X(\infty)=
		\begin{cases}
		\frac{1}{\Psi_X'(0+)}~~~~~~~~~~\bP\rbra{\lim_{t \uparrow \infty} 
		{X_t}=\infty}=1\\
		\infty~~~~~~~~~~~~~~~\text{otherwise}
		\end{cases}
		\end{align}
	(see e.g., \cite[pp.247]{Kyp}), we have
		\begin{align}
		\eqref{328a}
		%=&c\rbra{\Psi_X'(0+) \lor 0}
		%+c\int_0^\infty dv \int_{(-\infty, -v)}
		%	\rbra{e^{-\Phi_{X}(0)v} - e^{\Phi_{X}(q)u}}\Pi_X(du)\\
		=&c\rbra{\Psi_X'(0+) \lor 0}
		+c\int_{(-\infty, 0)}\Pi_X(du)\int_0^{-u}
			\rbra{e^{-\Phi_{X}(0)v} - e^{\Phi_{X}(q)u}}dv.
		\label{331a}
		\end{align}
%	where in \eqref{331a} we used the Fubini's theorem.\\
	We divide the remainder of the proof into two parts. 
	\par
	$(i)$ Suppose $\Psi_X'(0+)>0$. In this case, we have $\Phi_X(0)=0$ and so 
			\begin{align}
			\eqref{331a}
			%=&c\Psi_X'(0+)
			%+c\int_{(-\infty, 0)}\Pi_X(du)\int_0^{-u}\rbra{1 - e^{\Phi_{X}(q)u}}dv\\
			=&c\Psi_X'(0+)
			+c\int_{(-\infty, 0)}\rbra{ue^{\Phi_{X}(q)u} - u}\Pi_X(du)	.
			\label{336a}
			\end{align}
	Using \eqref{derivative of Laplace exponent}, we have
			\begin{align}
			\eqref{336a}%=&c\rbra{\gamma_X +\int_{(-\infty, -1]}u\Pi_X(du)}
				%+c\int_{(-\infty, 0)}\rbra{ue^{\Phi_{X}(q)u} - u}\Pi_X(du)	\\
			=&c\rbra{\gamma_X + \int_{(-\infty, 0)}
				\rbra{ue^{\Phi_{X}(q)u} - u1_{(-1, 0)}(u)}\Pi_X(du)}.
			\end{align}
	Using \eqref{326b}, we obtain $c=1$. \par
	$(ii)$ Suppose $\Psi_X'(0+)\leq0$. In this case, we have
			\begin{align}
			\eqref{331a}
			%=&c\int_{(-\infty, 0)}\Pi_X(du)\int_0^{-u}
			%	\rbra{e^{-\Phi_{X}(0)v} - e^{\Phi_{X}(q)u}}dv\\
			&=c\int_{(-\infty, 0)}\rbra{ue^{\Phi_{X}(q)u}+\frac{1}{\Phi_X (0)}
				-\frac{1}{\Phi_X (0)}e^{\Phi_X (0)u}}\Pi_X(du).
			\label{340a}
			\end{align}
		Since $\Psi(\Phi(0))=0$ and by \eqref{Laplace exponent} with $\sigma_X=0$, 
	we have
		\begin{align}
		\eqref{340a}=c\rbra{\gamma_X + \int_{(-\infty, 0)}
				\rbra{ue^{\Phi_{X}(q)u} - u1_{(-1, 0)}(u)}\Pi_X(du)}.
		\end{align}
	Using \eqref{326b}, we obtain $c=1$. 
	Thus the proof is complete.
	}
	We need the following two lemmas for later use. 
	\begin{Lem}\label{excursion Gerber-Shiu under a}
	For all $a > 0$, $q \geq 0$ and non-negative measurable function $f$ we have
		\begin{align}
		n^{X}\rbra{e^{-q\tau^{-}_{0}}f\rbra{X_{\tau^{-}_{0}}, X_{\tau^{-}_{0}-}}
			; \tau^{-}_{0}< \tau^{+}_{a}}
		=\int f(u, v)
			\overline{K}_X^{(q, a)} (du~ dv), 
		\label{345p}
		\end{align}
	where $\overline{K}_X^{(q, a)}$ is the measure carried on $(-\infty, 0)\times (0, \infty)$ defined by
		\begin{align}
		\overline{K}_X^{(q,a)}(du~ dv) =\overline{K}_X^{(q)}(du~ dv)
		:= \frac{W_{X}^{(q)}(a-v)}{W_{X}^{(q)}(a)}\tilde{\Pi}_{X}(du~ dv)
		\end{align}
	\end{Lem}
	The proof is parallel to that of {$(i)$ of Lemma \ref{Gerber-Shiu}}, so that we omit it. 
	\begin{Lem}\label{excursion potential density}
	For all $q \geq 0$ and non-negative measurable function $f$, we have
		\begin{align}
		n^{X}\rbra{ \int_0^{\tau^{-}_{0}\land T_0} e^{-qt}f\rbra{X_t}dt}
		=\int_0^\infty e^{-\Phi_{X}(q)y}f(y)dy. 
		\label{356p}
		\end{align}
	\end{Lem}
	\Proof{
	Using the monotone convergence theorem, we have
		\begin{align}
		n^{X}\rbra{ \int_0^{ \tau^{-}_{0} \land T_0} e^{-qt}f(X_t)dt}
		=&\lim_{\epsilon \downarrow 0}
			n^{X}\rbra{ \int_{\tau^{+}_{\epsilon}}^{ \tau^{-}_{0} \land T_0} e^{-qt}f(X_t)dt ; 
			\tau^+_\epsilon < \infty}
		\label{323}
		\end{align}
		and using the strong Markov property, we have
		\begin{align}
		\eqref{323}
		=&\lim_{\epsilon \downarrow 0}
			n^{X}\rbra{ e^{-q\tau^{+}_{\epsilon}}
			\rbra{%1_{\{X_{\tau^{-}_{0}-}> \epsilon\}}
			\int_{0}^{ \tau^{-}_{0} \land T_0} e^{-qt}f(X_t)dt
			}\circ \theta_{\tau^{+}_{\epsilon}}
			 ; \tau^+_\epsilon < \infty}
		\\
		=&\lim_{\epsilon \downarrow 0}
			n^{X}\rbra{ e^{-q\tau^{+}_{\epsilon}}; \tau^+_\epsilon < \infty }
			\bE^{X}_{\epsilon}\rbra{ \int_{0}^{ \tau^{-}_{0} \land T_0} e^{-qt}f\rbra{X_t}dt }
		\\
		=&\lim_{\epsilon \downarrow 0}
			\frac{1}{W_X^{(q)}(\epsilon)}
			\bE^{X}_{\epsilon}\rbra{ \int_{0}^{ \tau^{-}_{0} \land T_0} e^{-qt}f\rbra{X_t}dt },\label{3026}
		\end{align}
	where in \eqref{3026} we used Theorem \ref{equation lemma}.
	Using \eqref{one-sided potential density} with $b=0$, we obtain
		\begin{align}
		\int_{\epsilon}^{ \infty}f( y)e^{-\Phi_{X}(q)y} dy \leq
			\frac{1}{W^{(q)}_{X}
			(\epsilon)}\bE^{X}_{\epsilon}\rbra{ \int_{0}^{ \tau^{-}_{0}} e^{-qt}f\rbra{X_t}dt }
		\leq& \int_{0}^{ \infty}f( y)e^{-\Phi_{X}(q)y} dy.
		\end{align}
%	so we have 
%		\begin{align}
%		\lim_{\epsilon \downarrow 0}\int_{(\epsilon, \infty)}
%			f( y)e^{-\Phi_{X}(q)y}dy
%		\leq&\lim_{\epsilon \downarrow 0}
%			\frac{1}{W_X^{(q)}(\epsilon)}
%			\bE^{X}_{\epsilon}\rbra{ \int_{(0, \tau^{-}_{0})} e^{-qt}f\rbra{X_t}dt }
%		\\
%		\leq&\int_{(0, \infty)}
%			f( y)e^{-\Phi_{X}(q)y} dy.
%		\end{align}
		By the monotone convergence theorem, 
		the proof is complete.
	}
%	\begin{Rem}
%	The formulae \eqref{303p}, \eqref{317p}, \eqref{345p} and \eqref{356p}
%	can be found also in \cite{Avr}, \cite{Par2}, \cite{Par2} and \cite{Par}, 
%	respectively. 
%	We suspect their proofs contain ambiguity of a multiplicative constant of $n^X$. 
%	We believe our proofs contain no such ambiguity because we fixed 
%	the normalization of $n^X$ by \eqref{normalization of excursion}. 
%	\end{Rem}
	{\bf{(II)}} We assume that $X$ has bounded variation paths. 
	Note that in this case $0$ is irregular for itself. We write 
	\begin{align}
	n^X = \delta_X \bP^{X^0}_0.
	\end{align} 
	Then we have
		\begin{align}
		n^{X}\rbra{ e^{-q\tau_{a}^{+}}; \tau^+_a < \infty}=
		\delta_X \bE^{X^0}_0\rbra{ e^{-q\tau_{a}^{+}}; \tau^+_a < \infty}
		=\delta_X\frac{W_X^{(q)}(0)}{W_X^{(q)}(a)}=\frac{1}{W_X^{(q)}(a)}, 
		\label{526e}
		\end{align}
	where we used \eqref{exit problem} and \eqref{scale function 0}. 
	Thus we see that Theorem \ref{equation lemma} still holds in this case. 
	Lemmas \ref{excursion Gerber-Shiu}, 
	\ref{excursion Gerber-Shiu under a} and 
	\ref{excursion potential density}
	still hold as they are by a similar argument. 
	In particular, we obtain
	\begin{align}
	n^X\rbra{1-e^{-qT_0}}=\delta_X\bE^{X^0}_0\rbra{1-e^{-qT_0}}
	=\frac{1}{\Phi_X'(q)}={\Psi_X'}(\Phi_X(q)).
	\end{align}
	which may be regarded as the counterpart of 
	the normalization \eqref{normalization of excursion} 
	in the unbounded case.

\section{Kyprianou--Loeffen's Refracted L\'{e}vy processes}

Let us recall some results of Kyprianou--Loeffen \cite{Kyp2}. 
We fix a constant $\alpha >0$ and let $X$ be a general spectrally negative L\'evy 
process, which may possibly have Gaussian component. 
Set $Y_t =X_t +\alpha t$. 
Note that $0<\delta_X<\delta_X + \alpha =\delta_Y$ 
if $X$ has bounded variation paths. 
	\begin{Thm}[\cite{Kyp2}]\label{Kyp ref exi}
	%Suppose Assumption \ref{ref ass} is satisfied.
	For a fixed starting point $U_{0} = x \in \bR$, there exists a unique strong solution 
	to \eqref{Kyprianou--Loeffen's refracted}． 
	\end{Thm}
Let $U$ be a solution to Kyprianou--Loeffen's stochastic differential equation 
\eqref{Kyprianou--Loeffen's refracted}. 
	\begin{Thm}[\cite{Kyp2}]
	For all $ x \in [b, a]$ and $q \geq 0$, we have \eqref{1.05q} 
where $W_U^{(q)}$ is defined by \eqref{1.06q}. 
	\end{Thm}
They also calculated the potential densities with and without barriers. 
	\begin{Thm}[\cite{Kyp2}]
	For all $x \in [b, a]$, $q > 0$, we have 
		\begin{align}
		&\overline{\underline{r}}^{(q)}_{U}( x, y) 
		= 
		\begin{cases}
		 \frac{W_U^{(q)}(x, b)}
		{W_U^{(q)}(a, b)} W_{X}^{(q)}(a - y)-
		W_{X}^{(q)}(x - y) ~~~~~~~&y \in (0, a] 
\\
		 \frac{W_U^{(q)}
		(x, b)}{W_U^{(q)}(a, b)}
		 W_U^{(q)}(a, y) - W_U^{(q)}(x, y)  ~~~~~~~~~~~~&y \in [b, 0], 
		\end{cases}
\end{align}
\begin{align}
\un{r}_U^{(q)}(x, y)=
\begin{cases}
\frac{W_U^{(q)} (x, b)}{\alpha \un{H}_U^{(q)}(b)}e^{-\Phi_X(q)y}-W_X^{(q)}(x-y)\quad &y\in(0, \infty)\\
\frac{\un{H}_U^{(q)}(y)}{\un{H}_U^{(q)}(b) }W_U^{(q)}(x, b)-W_U^{(q)}(x, y)\quad &y \in [b, 0]
\end{cases}
\n 
\end{align}
with $\un{H}_U^{(q)} (y)=\int_0^\infty e^{-\Phi_X(q)z}W_Y^{(q)\prime}(z-y)dz$, 
\begin{align}
\bar{r}_U^{(q)}(x, y)=
\begin{cases}
\frac{\bar{H}_U^{(q)} (x)}{\bar{H}_U^{(q)} (a)}W_X^{(q)}(a-y)-W_X^{(q)}(x- y)\quad & y\in (0, a]\\
\frac{\bar{H}_U^{(q)} (x)}{\bar{H}_U^{(q)} (a)}W_U^{(q)}(a,y)-W_U^{(q)}(x, y)\quad & y\in (-\infty , 0]
\end{cases}
\n
\end{align}
with $\bar{H}_U^{(q)} (x)=e^{\Phi_Y (q)x}+\alpha \Phi_Y(q) \int_0^x e^{\Phi_Y (q)z}W_X^{(q)} (x-z)dz$, and 
\begin{align}
r_U^{(q)}(x, y)=
\begin{cases}
\frac{1}{\alpha}H_U^{(q)}(x)e^{-\Phi_X (q)y}-W_X^{(q)}(x-y )\quad &y \in ( 0, \infty)\\
H_U^{(q)}(x)\un{H}_U^{(q)}(y)-W_U^{(q)}(x, y)\quad & y\in (-\infty, 0]
\end{cases}
\n
\end{align}
with $H_U^{(q)}(x)=\bar{H}_U^{(q)} (x) \frac{\Phi_X(q) -\Phi_Y(q)}{\Phi_Y(q)}$, where $W_U^{(q)}$ has been given in \eqref{1.06q}. 
	\end{Thm}

\section{Generalization of refracted L\'{e}vy processes}\label{General}

We  now generalize Kyprianou--Loeffen's refracted L\'evy processes. 
We assume that $X$ and $Y$ 
are spectrally negative L\'{e}vy processes.
We assume, in addition, that 
	\begin{align}
	\text{$X$ has no Gaussian component whenever $X$ has unbounded variation paths.}
	\end{align}

In the unbounded variation case, 
we define the law of the stopped process $\bP^{U^{0}}_x$ by
	\begin{align}
&\bP^{U^0}_x\rbra{ F\rbra{(U_t)_{t < \tau^{-}_{0}}, (U_{t+\tau^{-}_{0}})_{t\geq 0}}}
=\bP^{X}_x\rbra{ \mathbb{E}^{Y^0}_{y}
\rbra{F\rbra{w, (Y^0_t)_{t\geq 0}}}
{\biggr|}_{\tiny{\begin{subarray} xy=X(\tau^{-}_{0})\\w=(X(t))_{t<\tau^{-}_{0}}\end{subarray}}}}~~~x \neq 0 
%	\label{law of stopped process of U}
	\end{align}
and the excursion measure $n^U$ by 
	\begin{align}
&n^{U}\rbra{ F\rbra{(U_t)_{t < \tau^{-}_{0}}, (U_{t+\tau^{-}_{0}})_{t\geq 0}}}
=n^{X}\rbra{ \mathbb{E}^{Y^0}_{y}
\rbra{F\rbra{w, (Y^0_t)_{t\geq 0}}}
{\biggr|}_{\tiny{\begin{subarray}
xy=X(\tau^{-}_{0})\\w=(X(t))_{t<\tau^{-}_{0}}\end{subarray}}}}
%	\tag{\ref{excursion measure of U}}
	\end{align}
for all non-negative measurable functional $F$, 
where $Y^0_t=Y_{t\land T_0}$ denotes the stopped process of $Y$ upon hitting zero.
Thus, we appeal to the excursion theory {(see Section \ref{construction})}, 
to construct the strong Markov process $U$ 
without stagnancy at $0$ (that is, $R_U^{(1)}1_{\{0\}}=0$) 
from $n^U$ together with $\{\bP^{U^0}_x\}_{x\neq0}$.
\par
In the bounded variation case, 
we define $U$ as a solution of \eqref{Kyprianou--Loeffen's refracted} 
constructed connecting $X$ and $Y$ mutually 
(this argument is similar as \cite{Kyp2}).
When $X$ and $Y$ are compound Poisson processes, uniqueness of the
solution of \eqref{Kyprianou--Loeffen's refracted} is easily proved.  
We write 
\begin{align}
n^X=\delta_X\bP^{X^0}_0 \text{ and } n^U=\delta_X\bP^{U^0}_0.
\end{align}
Then we obtain \eqref{excursion measure of U} as a formula. 
Therefore we can do a simultaneous discussion in between the two cases 
of bounded and of unbounded variation.
%We will think about potential measures of $U$ in section. \ref{sec pot}.

%If $X$ has bounded paths, we assume that $n^U := \bE^{U}_{0}$.
	\begin{Thm}\label{ref potential density}
	For all $q > 0$ and non-negative measurable function $f$ with $f(0)=0$, 
	we have
		\begin{align}
		N_U^{(q)}f:=&
		n^U\rbra{ \int_0^{T_0} e^{-qt}f(X_t)dt}\\
		=&\int_{0}^{ \infty}e^{-\Phi_{X}(q)y}f(y)dy
		+\int R_{Y^{0}}^{(q)}f(u)K^{(q)}_{X}(du~ dv).
		\label{506v}
		\end{align}
		%For all measurable function $f$ on $\bR$, we define $f(\Delta) = 0$.
	Consequently we have
		\begin{align}
		R_U^{(q)}f(0)&=\frac{N_U^{(q)}f}{qN_U^{(q)}1}, 
		\label{507v}\\
		R^{(q)}_{U}f(x)&=R_{Y^{0}}^{(q)}f(x)+e^{\Phi_Y(q)x}R_{U}^{(q)}f(0), ~~x <0,
		\label{504h}  \intertext{and}
		R^{(q)}_{U}f(x)&=\underline{R}_X^{(q; 0)}f(x)
			+\int R_{U}^{(q)}f(u)G_X^{(q)}(x, du~ dv), ~~x>0,
		\label{505h}
		\end{align}
		where
		\begin{align}
		\underline{R}_X^{(q; 0)}f(x)
		=\bE^X_x\rbra{\int_0^{\tau^-_0} e^{-qt} f(X_t)dt}.
%		R^{Y^0}_{q}f\left(u\right) = \mathbb{E}^{Y}_{u}\rbra{
%		\int_{(0, T_0)}e^{-qt}f(Y_t)dt }. %\\
%		K^{X}_{q}g=\int_{(-\infty, 0)\times (0, \infty)}g(u)K^{X}_{q}(du, dv).
		\end{align}
	\end{Thm}
	\Proof{ 
	Let us calculate $N_U^{(q)}f$. 
	Since $\int_0^{T_0}=\int_0^{\tau_0^-}+\int_{\tau_0^-}^{T_0}$, we have that 
	$N_U^{(q)}f$ is equal to 
		\begin{align}
		n^{U}\rbra{\int_{0}^{\tau^{-}_0}e^{-qt}f(U_t)dt }
		+n^{U}\rbra{e^{-q\tau^{-}_{0}}\rbra{\int_{0}^{ T_0}e^{-qt}f(U_t)dt}\circ
			\theta_{\tau^{-}_{0}} ;\tau^{-}_0<\infty}.
		\label{422d}
		\end{align}
	Using Lemmas \ref{excursion potential density} and \ref{excursion Gerber-Shiu}, 
	we obtain \eqref{506v}. 
	\par
	Let us prove \eqref{507v}. 
	{Note that the finiteness of $N_U^{(q)}1$ will be proved in Lemma \ref{Lem304h}.} 
	When $X$ has unbounded variation paths, the formula \eqref{507v} can be found, 
	e.g., in \cite[pp.423]{Rog}. 
	Suppose $X$ has bounded variation paths. 
	We denote 
		$T_0^{(0)} = 0$ and define 
		\begin{align}
		T_0^{(n)}=\inf\cbra{t> T_0^{(n-1)}: X_t = 0} \nonumber
		\end{align} 
		recursively for 
		all $n \in \bN$. 
		Then we have
		\begin{align}
		R_U^{(q)}f(0)
		=&\sum_{n=0}^{\infty} \bE_0^{U}\rbra{ \int_{T_0^{(n)}}^{T_0^{(n+1)}}
			e^{-qt}f(U_t)dt; T_0^{(n)} < \infty}
		\\
		=&\sum_{n=0}^{\infty}\bE_0^{U}\rbra{e^{-qT_0}}^n
			\bE_0^U \rbra{ \int_0^{T_0} e^{-qt}f(U_t)dt}
		\\
		=&\frac{\bE_0^{U^0} \rbra{ \int_0^{T_0} e^{-qt}f(U_t)dt}}
			{q\bE_0^{U^0}\rbra{\int_0^{T_0}e^{-qt} dt}}.
		\label{525i}
		\end{align}
	Since we write $n^U = \delta_X\bP^{U^0}_0$, we obtain \eqref{507v}. 
	\par
	The remainder of the proof is straightforward. 
	}
	The following theorem shows the choice of $n^U$ leads to a normalization 
	similar to \eqref{normalization of excursion}. 
	\begin{Thm}
	For all $q>0$, we have 
		\begin{align}
	n^U\rbra{1 - e^{-qT_0}}&=\frac{1}{r_U^{(q)}(0, 0+ ) } \label{normalization}\\
	&=\rbra{ \Psi_{X}^{'}(0) \lor 0}
			+\int
			\rbra{e^{\Phi_{X}(0)u} - e^{\Phi_{Y}(q)u-\Phi_{X}(q)v}}
			{\tilde{\Pi}}_{X}(du ~dv).
		\label{515pika}
	\end{align}
	\end{Thm}
	\Proof{
	By \eqref{507v} of Theorem \ref{ref potential density}, we have 
	\begin{align}
	r_U^{(q)} (0, y) =\frac{1}{q N_U^{(q)} 1}
			\rbra{e^{-\Phi_X(q)y} 1_{(y>0)} + \int r_{Y^0}^{(q)} (u, y)K_X^{(q)}(du~dv)}. 
	\end{align}
	Since $r_{Y^0}^{(q)}(u, y)=0$ for ${u}<0$ and $y>0$, we have 
	\begin{align}
	r_U^{(q)}(0, 0+)=\frac{1}{qN_U^{(q)}1}. 
	\end{align}
	On the other hand, we have
	\begin{align}
	q N_U^{(q)} 1 =n^U\rbra{1 - e^{-qT_0}} 
	\label{518aa}
	\end{align}
	by the definition of $N_U^{(q)}$.  
	Thus we obtain \eqref{normalization}. \par 
	The other expression 
	\eqref{515pika} can be proved easily by a similar argument to \eqref{331a}. 
	}

%	\begin{Rem}
%	$n^U\rbra{1-e^{-qT_0} }$ can be represented as \eqref{518aa} using 
%	scale functions and Laplace exponents of $X$ and $Y$. 
%	And $n^U\rbra{1-e^{-qT_0} }$ has another representation by making a computation %similar to Lemma \ref{constant}: 
%			\begin{align}
%		&n^{U}\rbra{1-e^{-q T_0}}\nonumber\\
%		=&n^{X}\rbra{\tau^{-}_{0}= \infty}
%		+n^{X}\rbra{1- e^{-q\tau^{-}_{0}}\bE^{Y}_{X(\tau^{-}_{0})}
%			\rbra{e^{-qT_0}};\tau^-_0 < \infty}
%		\\
%		=&\rbra{ \Psi_{X}^{'}(0) \lor 0}
%			+\int_{0}^{\infty}dv\int_{(-\infty, 0)}
%			\rbra{e^{\Phi_{Y}(0)u} - e^{\Phi_{Y}(q)u-\Phi_{X}(q)v}}
%			\Pi_{X}(du - v).
%		\label{436g}
% 		\end{align} 
%	\end{Rem}

\section{Exit Problem of generalized refracted L\'{e}vy processes}

We prepare a general formula.
	\begin{Lem}\label{excursion exit}
	Let $Z$ be a standard process
	with no positive jumps without stagnancy at $0$ (i.e., $R_Z^{(1)}1_{\{0\}}=0$). 
	If $0$ is regular for itself, then
		\begin{align}
		\bE^{Z}_{0}\rbra{ e^{-q\tau_{a}^{+}} ;\tau_{a}^{+} <\tau_{b}^{-}}
		=\frac{n^{Z}\rbra{ e^{-q\tau_{a}^{+}}; \tau^+_a < \infty }}
			{n^{Z}\rbra{ 1 - e^{-qT_0}1_{\{ \tau^+_a = \infty , \tau^-_b = \infty \}}}}
		\label{601j}
		\end{align}
	for all $a > 0 > b$ and $q \geq 0$, where $n^Z$ denotes an excursion measure 
	away from $0$.
	If $0$ is irregular for itself, the identity \eqref{601j} still holds where $n^Z$ 
denotes a constant multiple of $\bP^{Z^0}_0$.
	\end{Lem}
	\Proof{
	{It is sufficient to prove \eqref{601j} only when $q>0$.} 
	We assume first that $0$ is regular for itself. 
	Let $p$ denote a Poisson point process 
	{defined on the probability space $(\Omega , {\cal{F}} , \bP)$} 
	with characteristic measure $n^Z$.  
	Set $\eta(s)=\sum_{u\leq s}T_0\rbra{p(u)}$. 
	{Note that $\eta$ will be the inverse local time at $0$ for the process 
	constructed from the excursions, which equals in law to $Z$ under $\bP^Z_0$.}
	For $E\in {\cal{B}} (\bD)$, we write $\kappa_E =\inf\{ s \geq 0: p(s)\in E\}$. 
	We let {$A=\{\tau^+_a < \infty\}\cup\{\tau^-_b < \infty\}$} 
	and we denote by $\epsilon^{\ast}=p(\kappa_{A})$ 
	the first excursion belonging to $A$. 
	Then we have   
		\begin{align}
		\bE^{Z}_{0}\rbra{ e^{-q\tau_{a}^{+}} ;\tau_{a}^{+} <\tau_{b}^{-}}
		=&\bE\rbra{ e^{-q\eta(\kappa_{A }-)}
			e^{-q\tau^+_a\rbra{\epsilon^\ast}}}\\
		=&\bE\rbra{ e^{-q\eta(\kappa_{A }-)}}
			\frac{n^{Z}\rbra{ e^{-q\tau_{a}^{+}} ;A}}
			{n^{Z}(A)}
			\label{q>0} \\
		=&{\bE\rbra{ e^{-q\eta(\kappa_{A }-)}}
			\frac{n^{Z}\rbra{ e^{-q\tau_{a}^{+}} ;\tau^+_a < \infty }}{n^{Z}(A)} }\label{604ac}
		\end{align}
	where {$\bE$ denotes the expectation with {respect} to $\bP$. 
	Note that in \eqref{q>0}} 
	we used the renewal property of the Poisson point process 
	{and in \eqref{604ac} we used the fact that 
	$e^{-q\tau^+_a}=0$ on $\{\tau^+_a=\infty,  \tau^-_b < \infty \}$. 
%	the same argument as Remark \ref{Rem303}.
	} 
	We write {$p_{A^c} $} 
	for $p$ restricted to excursions belonging  to ${A}^c$ 
	and write $\eta_{A^c}(s)=\sum_{u\leq s}T_0 \rbra{p_{A^c}(u)}$. 
	Since $\eta (\kappa_{A}-)=\eta_{A^c} (\kappa_{A})$ where $\eta_{A^c}$ and 
	$\kappa_{A}$ are independent, we have 
		\begin{align}
		\bE\rbra{ e^{-q\eta(\kappa_{A}-)}}
		=&n^{Z}(A)\int_{0}^{ \infty}e^{-n^Z(A)t}
			{\bE} \rbra{ e^{-q\eta_{A^c}(t)} }dt
		\\
		=&n^{Z}(A)\int_{0}^{ \infty}e^{-n^Z(A)t}
			\rbra{\exp(-tn^{Z}(1-e^{-qT_0}; {A}^c))}dt
		\\
		=&\frac{n^{Z}\rbra{ A}}
			{n^{Z}\rbra{ 1 - e^{-qT_0}1_{A^{c}}}}.
			\label{q kappa}
		\end{align}
	Thus we obtain \eqref{601j}. 
	\par
	We second assume that $0$ is irregular for itself.
	Using the notation of the proof of Theorem \ref{ref potential density}, we have
		\begin{align}
		\bE_0^{Z}\rbra{e^{-q\tau^+_a} ; \tau^+_a < \tau^-_b}
		=&\sum_{n=0}^{\infty}
			\bE_0^{Z}\rbra{e^{-q\tau^+_a} ; T^{(n)}_0 <\tau^+_a < \rbra{T^{(n+1)}_0\land\tau^-_b}}\\
		=&\sum_{n=0}^{\infty}
			{\bE_0^{Z^0}}\rbra{e^{-qT_0}; \tau^+_a = \infty , \tau^-_b = \infty}^n
			\bE_0^{Z^0}\rbra{e^{-q\tau^+_a} ; \tau^+_a < \infty}\\
		=&\frac{\bE_0^{Z^0}\rbra{e^{-q\tau^+_a} ; \tau^+_a < \infty}}
			{1-\bE_0^{Z^0}\rbra{e^{-qT_0}; \tau^+_a = \infty , \tau^-_b = \infty}}.
		\label{613i}
		\end{align}
	Thus we obtain \eqref{601j}.
	}
	\begin{Thm}
	\label{refracted exit problem}
	For all $x \in [b, a]$ and $q \geq 0$, we have
		\begin{align}
		\mathbb{E}^{U}_{x}\rbra{e^{-q\tau^{+}_{a}} ; \tau^{+}_a < \tau^{-}_b }
		=\frac{W_{U}^{(q)}(  x, b )}{W_{U}^{(q)}(  a, b )} , 
		\label{610z}
		\end{align}
	where {the function $W_U^{(q)}(x, y)$ is defined as follows: for $x \in (0, \infty)$, 
	\begin{align}
	W_U^{(q)}(  x, y )=&
	W^{(q)}_{X}(x)W^{(q)}_{Y}(-y)\rbra{\Psi_{X}^{'}(0) \lor 0 }\\
		&+\int
		{\big(} W^{(q)}_{X}(x)W^{(q)}_{Y}(-y)e^{\Phi_{X}(0)u} 
		- W^{(q)}_{Y}(u-y)W^{(q)}_{X}(x-v){\big)}\tilde{\Pi}_{X}(du~ dv)
	\end{align}
	and for $x \in (-\infty , 0]$, 
	\begin{align}
	W_U^{(q)}(  x, y )=W_{Y}^{(q)}(x-y). 
	\end{align}}
	\end{Thm}
	\Proof{
	We discuss the two cases of unbounded and of bounded variation at the same time. 
	\par
	We calculate 
	$\mathbb{E}^{U}_{0}\rbra{e^{-q\tau^{+}_{a}} ; \tau^{+}_a < \tau^{-}_b }$.
	Using Lemma \ref{excursion exit}, we have
		\begin{align}
		&\mathbb{E}^{U}_{0}\rbra{ e^{-q\tau_{a}^{+}}
 ;\tau_{a}^{+} <\tau_{b}^{-}}
		=\frac{n^{U}\rbra{ e^{-q\tau_{a}^{+}}; \tau^+_a < \infty}}
			{n^{U}\rbra{ 1 - e^{-qT_0}
			1_{\{ \tau^+_a = \infty , \tau^-_b = \infty \}}}}.
		\end{align}
		Using Theorem \ref{equation lemma}, we can rewrite the numerator as 
		\begin{align}
		n^{U}\rbra{ e^{-q\tau_{a}^{+}}; \tau^+_a < \infty}=
		n^{X}\rbra{ e^{-q\tau_{a}^{+}}; \tau^+_a < \infty}=\frac{1}{W_{X}^{(q)}(a)}.
		\label{416}
		\end{align}
		We divide the denominator $n^{U}\rbra{ 1 - e^{-qT_0}
			1_{\{ \tau^+_a = \infty , \tau^-_b = \infty \}}}$ into the following sum: 
		\begin{align}
			{n^{U}\rbra{1-e^{-q T_0}}+
			n^{U}\rbra{ e^{-qT_0};\{ \tau^+_a < \infty\}\cap\{\tau^-_b = \infty \}}
			+n^{U}\rbra{ e^{-qT_0};{\tau^-_b < \infty }} }.
		\label{415}
		\end{align}
	Let us compute these expectations. 
%	When $X$ is irregular for $0$, $X$ has bounded variation paths and we have 
%		\begin{align}
%		n^{U}\rbra{ e^{-q\tau_{a}^{+}}; \tau^+_a < \infty}=
%		\delta_X \bE^{X}_0\rbra{ e^{-q\tau_{a}^{+}}; \tau^+_a < \infty}
%		=\delta_X\frac{W_X^{(q)}(0)}{W_X^{(q)}(a)}=\frac{1}{W_X^{(q)}(a)}, 
%		\label{526e}
%		\end{align}
%	where in \eqref{526e} we used \eqref{exit problem} and \eqref{scale function 0}.
	For the second term, we have  
		\begin{align}
		&n^{U}\rbra{ e^{-qT_0};\{ \tau^+_a < \infty\}\cap\{\tau^-_b = \infty \}}
		\nonumber
		\\
		=&n^{U}\rbra{ e^{-q\tau^{+}_{a}}\rbra{e^{-qT_0}1_{\{\tau^-_b = \infty \}}}
			\circ \theta_{\tau^{+}_{a}}; \tau^+_a < \infty}
		\\
		=&n^{X}\rbra{ e^{-q\tau_{a}^{+}}; \tau^+_a < \infty}
			\mathbb{E}^{X}_{a}\rbra{ e^{-q\tau^{-}_{0}}\mathbb{E}^{Y}_{X(\tau^{-}_{0})}
			\rbra{ e^{-qT_{0}}1_{\{T_0 < \tau^{-}_{b}\}}}
			;\tau^{-}_{0} < \infty}
			\label{419}
		\\
		=&%n^{X}\rbra{ e^{-q\tau_{a}^{+}}; \tau^+_a < \infty}
			\frac{1}{W_X^{(q)}(a)}
			%\int_{0}^{ \infty}dv
			\int
			\frac{W_{Y}^{(q)}(u-b)}{W_{Y}^{(q)}(-b)}G^{(q)}_X(a, du~ dv), 
%		\rbra{e^{-\Phi_{X}(q)v} W_{X}^{(q)}(a) -W_{X}^{(q)}(a - v)}\Pi_{X}(du - v)
		\label{421}
		\end{align}
	where in \eqref{421} we used Theorem \ref{equation lemma},
	{$(i)$ of Theorem \ref{Gerber-Shiu}} and \eqref{exit problem}. 
	For the third term, 
	we have
		\begin{align}
		n^{U}\rbra{ e^{-qT_0};\tau^-_b < \infty }
		=&n^{U}\rbra{e^{-q\tau^{-}_{0}}\mathbb{E}^{U}_{U(\tau^{-}_{0})}
			\rbra{ e^{-qT_0};\tau^{-}_{b}<T_0}	;\tau^{-}_{0} < \infty }. 
		\label{425}
	\end{align}
	Using Lemma \ref{excursion Gerber-Shiu}, we have
	\begin{align}
		\eqref{425}=&\int e^{-\Phi_{X}(q)v}\mathbb{E}^{Y}_{u}
			\rbra{ e^{-qT_0};\tau^{-}_{b}<T_0} \tilde{\Pi}_{X}(du ~ dv)
		\label{426}
		\\
		=&\int e^{-\Phi_{X}(q)v}
		\rbra{e^{\Phi_{Y}(q)u}-\frac{W_{Y}^{(q)}(u-b)}{W_{Y}^{(q)}(-b)}}
		\tilde{\Pi}_{X}(du ~ dv),
		\label{427}
	\end{align}
	where in \eqref{427} we used \eqref{exit problem} and \eqref{one-sided exit problem}.
	Therefore, using \eqref{515pika}, we obtain
%		\begin{align}
%		{\eqref{415}}=
%		\frac{W_U^{(q)}(  0, b )}{W_U^{(q)}( a, b )}.
%		\label{436}
%		\end{align}
	\begin{align}
	n^U\rbra{ 1- e^{-qT_0}1_{\{ \tau^+_a = \infty, \tau^-_b = \infty \}}}
	=\frac{1}{W_X^{(q)}(a)}\frac{W_U^{(q)}(  a, b )}{W_U^{(q)}( 0, b )}
	\label{436}
	\end{align}
	and we obtain \eqref{610z} for $x=0$.
	For all $x < 0$, we have
		\begin{align}
		\bE^{U}_{x}\rbra{e^{-q\tau^{+}_{a}} ; \tau^{+}_a < \tau^{-}_b }
		=&\bE^{Y}_{x}\rbra{ e^{-q\tau^{+}_{0}}; \tau^{+}_{0}<\tau^{-}_{b}}
			\bE^{U}_{0}\rbra{e^{-q\tau^{+}_{a}} ; \tau^{+}_a < \tau^{-}_b }.
		\label{438}
		\end{align}
	Using \eqref{exit problem} and \eqref{610z} for $x=0$, 
	we have \eqref{610z} for $x<0$. 
	For all $x>0$, we have
		\begin{align}
		&\bE^{U}_{x}\rbra{e^{-q\tau^{+}_{a}} ; \tau^{+}_a < \tau^{-}_b }\nonumber
		\\
		=&\bE^{U}_{x}\rbra{e^{-q\tau^{+}_{a}} ; \tau^{+}_a < \tau^{-}_0 }
		+\bE^{U}_{x}\rbra{e^{-q\tau^{-}_{0}} 
		\rbra{e^{-q\tau^{+}_{a}}1_{\{\tau^{+}_{a}<\tau^{-}_{b}\}}}
		\circ \theta_{\tau^{-}_{0}}; \tau^{-}_0 < \tau^{+}_a }\\
		=&\frac{W_{X}^{(q)}(x)}{W_{X}^{(q)}(a)}
		+\bE^{X}_{x}\rbra{e^{-q\tau^{-}_{0}} 
		\bE^{Y}_{X(\tau^{-}_{0})}
		\rbra{e^{-q\tau^{+}_{a}};\tau^{+}_{a}<\tau^{-}_{b}}
		; \tau^{-}_0 < \tau^{+}_a },
		\label{444}
		\end{align}
	where 
	in \eqref{444} we used \eqref{exit problem}.
	Using \eqref{610z} for $x<0$ and {$(ii)$ of Theorem \ref{Gerber-Shiu}}, 
	the second term of \eqref{444} is equal to 
		\begin{align}
		\int
		\frac{W_U^{(q)}( u, b)}{W_U^{(q)}( a, b)}
		\overline{\underline{r}}_X^{(q; 0, a)}(x, v)
		\tilde{\Pi}_{X}(du~ dv). 
		\label{452}
		\end{align}
	Thus we obtain \eqref{610z} for $x>0$. 
	The proof is complete. 
	}
	\begin{Cor}\label{refracted one-sided exit problem}
	For all $x \in (-\infty, a]$ and $q \geq 0$, we have 
		\begin{align}
		\bE^U_x\rbra{ e^{-q\tau^+_a}}
		=\frac{\overline{W}_U^{(q)} ( x)}{\overline{W}_U^{(q)} ( a)}
		\label{629z}
		\end{align}
	where {the function $\bar{W}_U^{(q)}(x)$ is defined as follows: 
for $x\in(0,\infty)$, 
		\begin{align}
		\overline{W}_U^{(q)} ( x)=&
		W^{(q)}_{X}(x)\rbra{\Psi_{X}^{'}(0) \lor 0 }\\
		&+\int
		{\big(} W^{(q)}_{X}(x)e^{\Phi_{X}(0)u}
		- W^{(q)}_{X}(x-v)e^{\Phi_Y(q)u}{\big)}\tilde{\Pi}_{X}(du~ dv)
		\end{align}
	and for $x \in (-\infty , 0]$, 
	\begin{align}
	\overline{W}_U^{(q)} ( x)=&e^{\Phi_Y(q)x}.
		\label{648k}
	\end{align}}
	In particular, 
	$\overline{W}_U^{(q)}( x)$ is a continuous and increasing function of $x$.
	\end{Cor}
	\Proof{
	Using the monotone convergence theorem and 
	Theorem \ref{refracted exit problem}, we have 
	\begin{align}
	\bE^U_x\rbra{ e^{-q\tau^+_a}}
	=\lim_{b \downarrow -\infty} \bE^U_x\rbra{e^{-q\tau^+_a} ;\tau^+_a<\tau^-_b}
	=\lim_{b \downarrow -\infty} 
		\frac{W_U^{(q)} (x, b)/ W_Y^{(q)}(-b) }{W_U^{(q)} (a, b)/ W_Y^{(q)}(-b)}.
	\end{align}
	Using the last equality of \cite[pp.124]{matter}, we have 
	\begin{align}
	\lim_{b \downarrow -\infty}\rbra{W_U^{(q)} ( x, b)/ W_Y^{(q)}(-b)}
	=\overline{W}_U^{(q)}( x), 
	\end{align}
	and we have \eqref{629z}. \par
	Next, we prove that $\overline{W}_U^{(q)}$ is 
	increasing and continuous. 
	It is obvious that $\overline{W}_U^{(q)}$ is increasing and continuous 
	on $(-\infty, 0]$, since $\overline{W}_U^{(q)}(x)=e^{\Phi_Y(q)x}$. 
	Using the dominated convergence theorem, we have 
	\begin{align}
	\lim_{\epsilon \downarrow 0}\overline{W}_U^{(q)}(\epsilon)
	=\lim_{\epsilon \downarrow 0}\frac{1}{\bE_0^U\rbra{e^{-q\tau^+_\epsilon}}}
	=\frac{1}{\bE_0^U\rbra{\lim_{\epsilon \downarrow 0}e^{-q\tau^+_\epsilon}}}=1,
	\end{align}
	so that we see $\overline{W}_U^{(q)}$ is continuous at $0$.
	Since  
	\begin{align}
	\overline{W}_U^{(q)}(x) =\frac{1}{\bE_0^U\rbra{e^{-q\tau^+_x}}},
	\end{align}
	it is thus sufficient to prove 
	that $\bE_0^U\rbra{e^{-q\tau^+_x}}$ is decreasing and continuous 
	on $(0, \infty)$. 
	For $0<x<y$, we have 
	\begin{align}
	\bE_0^U\rbra{e^{-q\tau^+_x}}-\bE_0^U\rbra{e^{-q\tau^+_y}}
	=&\bE_0^U\rbra{e^{-q\tau^+_x}}
		\rbra{1-\bE_x^U\rbra{e^{-q\tau^+_y}}}\geq0.
	\end{align}
	Using \eqref{exit problem}, for $x>0$, we have 
	\begin{align}
	\limsup_{\epsilon \downarrow 0}
	\absol{\bE_0^U\rbra{e^{-q\tau^+_{x-\epsilon}}}-
		\bE_0^U\rbra{e^{-q\tau^+_{x+\epsilon}}}}
	=&\limsup_{\epsilon \downarrow 0}\bE_0^U\rbra{e^{-q\tau^+_{x-\epsilon}}}\rbra{1-
				\bE_{x-\epsilon}^U\rbra{e^{-q\tau^+_{x+\epsilon}}}}\\
	\leq&\limsup_{\epsilon \downarrow 0}%\bE_0^U\rbra{e^{-q\tau^+_x}}
		\rbra{1-
			\bE_{x-\epsilon}^X\rbra{e^{-q\tau^+_{x+\epsilon}};\tau^+_{x+\epsilon}<\tau^-_0}}\\
	=&%\bE_0^U\rbra{e^{-q\tau^+_x}}
		\rbra{1-\lim_{\epsilon \downarrow 0}
				\frac{W_X^{(q)}(x-\epsilon)}{W_X^{(q)}(x+\epsilon)}}=0. 
	\end{align}
	The proof is complete.
	}

Let $C_0$ denote the set of continuous functions 
$ f : \bR \rightarrow \bR$ which vanish at $+\infty$ and $-\infty$. 
Note that $C_0$ is a Banach space with respect to the supremum norm 
$\norm{f} = \sup_{x\in\bR}\absol{f(x)}$ for $f\in C_0$.
\begin{Thm}\label{Feller property}
	Our generalized refracted L\'evy process is a Feller process.
\end{Thm}
	\Proof{
	Since $R_U^{(q)}$ comes from transition operators, 
	it is sufficient to verify the following conditions:\par
	$(i)$ For all $q>0$, $R_U^{(q)}$ is a map from $C_0$ to $C_0$. \par
	$(ii)$ For all $f\in C_0$, $\lim_{q\uparrow \infty} \norm{qR_U^{(q)}f-f}=0$.
	\par
	{\bf $1$) The proof of $(i)$}
	\par
	First, we prove that $R_U^{(q)}f$ is continuous. 
	Let $x\in\bR$ and $\epsilon >0$.
	Noting that $U$ has no positive jump, we have 
		\begin{align}
	&\absol{R_{U}^{(q)} f(x + \epsilon) - R_{U}^{(q)} f(x)}\nonumber \\
	&\leq 
	\absol{R_{U}^{(q)} f(x + \epsilon) 
		-\bE^{U}_x\rbra{e^{-q\tau^+_{x + \epsilon}}}R_{U}^{(q)} f(x + \epsilon)}
		+\absol{\bE^{U}_x \rbra{\int_0^{\tau^+_{x + \epsilon}}e^{-qt}f(U_t )dt}}\\
	&\leq
	\absol{R_{U}^{(q)} f(x + \epsilon)}
		\rbra{1- \bE^{U}_x\rbra{e^{-q\tau^+_{x + \epsilon}}}}
		+\norm{f}\bE^{U}_x \rbra{\int_0^{\tau^+_{x + \epsilon}} e^{-qt}dt} \\
	&\leq\frac{2}{q}\norm{f}\rbra{1-\bE^{U}_x\rbra{e^{-q\tau^+_{x + \epsilon}}}}. 
	\label{640pika}
	\end{align}
		By Corollary \ref{refracted one-sided exit problem}, we have
	\begin{align}
	\eqref{640pika}
	=\frac{2}{q}\norm{f}\rbra{1- \frac{\overline{W}_{U}^{(q)}( x)}
				{\overline{W}_{U}^{(q)}( x + \epsilon)}} 
	\rightarrow 0 \text{ as } \epsilon \downarrow 0.
	\label{641pika}
	\end{align}
	Thus we obtain right-continuity of $R_U^{(q)}f$.
	For the left-continuity {we have 
	\begin{align}
	&\absol{R_{U}^{(q)} f(x - \epsilon) - R_{U}^{(q)} f(x)} \\
	&\leq \absol{\bE^U_{x-\epsilon}\rbra{e^{-q\tau^+_x}}R_{U}^{(q)} f(x) - R_{U}^{(q)} f(x)}
	+\absol{\bE^{U}_{x-\epsilon} \rbra{\int_0^{\tau^+_{x }}e^{-qt}f(U_t )dt}}
	\end{align}
	and the remainder of its proof is similar to that of the right-continuity. }
	\par
	Second, we prove that $R_U^{(q)}f$ vanishes at $-\infty$. 
	For $x<0$, we may rewrite \eqref{504h} as 
	\begin{align}
	R_U^{(q)}f(x)=&R_{Y}^{(q)}f(x)-e^{\Phi_Y (q)x} R_{Y}^{(q)}f(0)
		+e^{\Phi_Y (q)x} R_{U}^{(q)}f(0).
	\label{536n}
	\end{align}
	By the Feller property of $Y$, we see that  
	$\lim_{x \downarrow -\infty}R_{U}^{(q)}f(x)=0$. 
	\par
	Third, we prove that $R_U^{(q)}f$ vanishes at $+\infty$. 
	We may assume without loss of generality that $f \geq 0$.
	For all $x > 0$, we have 
	\begin{align}
		R_U^{(q)}f(x)
		&=\bE^U_x\rbra{ \rbra{\int_0^{\tau^-_0} + \int_{\tau^-_0}^\infty} 
				e^{-qt}f(U_t)dt}\\
		&\leq R_X^{(q)}f(x) + \bE^X_x \rbra{e^{-q\tau^-_0}R_U^{(q)}f(X_{\tau^-_0})}
		\cdot \frac{1}{q}\norm{f}.
	\end{align}
	By the Feller property of $X$ and by the fact that 
	$\bE^X_x \rbra{e^{-q\tau^-_{0}}} = \bE^X_0 \rbra{e^{-q\tau^-_{-x}}}
	\rightarrow 0$ as $x \rightarrow \infty$, we obtain 
	$\lim_{x \uparrow \infty}R_U^{(q)}f(x) =0$.
	\par 
	{\bf $2$) The proof of $(ii)$}
	\par
	Define 
	\begin{align}
	{\omega}_\epsilon(f; x) = \sup_{y: \absol{y-x}\leq \epsilon}\absol{f(y) - f(x)}.
	\end{align} 
	Let us prove the pointwise convergence: 
		\begin{align}
		\lim_{q\uparrow \infty} qR_U^{(q)}f (x) = f(x), ~~~~~~x\in \bR.
		\label{642pikatyu}
		\end{align}
	For $x \in \bR$ and $\epsilon > 0$, we have 
	\begin{align}
	&\absol{ qR_U^{(q)}f (x) - f(x)}\\
	&\leq
	q\bE^U_x\rbra{\rbra{\int_0^{\tau^+_{x+\epsilon} \land \tau^-_{x-\epsilon}}
			+\int_{\tau^+_{x+\epsilon} \land \tau^-_{x-\epsilon}}^\infty}
				e^{-qt}\absol{f(U_t)-f(x)}dt} 
	\\
	&\leq
		\bE^U_x\rbra{1-e^{-q( \tau^+_{x+\epsilon} \land \tau^-_{x-\epsilon})}}
		{\omega}_{\epsilon}(f ;x) 
	+
		2\norm{f}\bE^U_x\rbra{e^{-q( \tau^+_{x+\epsilon} \land \tau^-_{x-\epsilon})}}.
		\label{645pikatyu}
	\end{align}
	We thus obtain 
	${\limsup}_{q \uparrow \infty} \absol{q R_U^{(q)} f(x) -f(x)}
		\leq {\omega}_{\epsilon}(f ; x)$
	for all $\epsilon >0$, 
	which proves \eqref{642pikatyu}. 
	By a standard argument with the help of the fact that the dual space of $C_0$ 
	can be identified with the space of signed measures, we can see that 
	$R_U^{(p)}(C_0)$ is dense in $C_0$ for all $p>0$. 
	\par
	Let $f = R_U^{(1)}g$ for some $g \in C_0$. 
	Using the resolvent equation, we have  
	\begin{align}
	\norm{f-q R_U^{(q)}f }
	=\norm{R_U^{(q)}g - R_U^{(q)}f}
	\leq \frac{1}{q}\norm{g - f} \rightarrow 0, \text{ as } q \uparrow \infty. 
	\label{647pikatyu}
	\end{align}
	Since $R_U^{(1)}(C_0)$ is dense in $C_0$, we obtain claim $(ii)$.  
	\par
	The proof is now complete. 
	}

\section{Potential measure of killed refracted L\'{e}vy processes}\label{sec pot}

In this section, we calculate the potential measure of refracted L\'{e}vy 
processes killed on exiting $[b, a]$.
	\begin{Thm}\label{killed ref potential density}
	For all $x \in [b, a]$ and $q \geq 0$, we have
		\begin{align}
		\overline{\underline{r}}^{(q)}_{U}( x, y)=
		\begin{cases}
		\frac{W_U^{(q)}(  x, b)}{W_U^{(q)}(  a, b)}W_X^{(q)}(a - y) - W_X^{(q)}(x - y),
			~~~~~~~&y \in (0, a]
		\\
		\frac{W_U^{(q)}( x, b)}{W_U^{(q)}( a, b)}W_U^{(q)}(  a, y)-W_U^{(q)}(  x, y),
			~~~~~~~&y \in [b, 0).
		\end{cases}
		\label{701z}
		\end{align}
	\end{Thm}
	\Proof{
	We follow the notation of Lemma \ref{excursion exit} for $L$, $\eta$, 
	$\kappa$, etc. \par
	{\bf{Step.1}} We calculate in the case $x = 0$.
	When $X$ has unbounded variation paths, we have
		\begin{align}
		&\bE^{U}_{0}\rbra{\int_{0}^{\tau^{+}_{a} \land\tau^{-}_{b}} e^{-qt}  f(U_{t} ) dt}
		\nonumber
		\\
		=&\bE^{U}_{0}\rbra{ 
		\int_{(0, \infty)}e^{-qs}1_{\{s < \eta(\kappa_{A} - )\}}dL(s)}
		n^{U}\rbra{\int_{0}^{T_0\land \tau^{+}_{a}\land\tau^{-}_{b}}e^{-qt}f(U_{t})dt
		}.
		\label{604}
		\end{align}
	where we used the compensation theorem of 
	the excursion point process.
	We may rewrite \eqref{604} using $\eta_{A^c}$, as	
	\begin{align}
		\bE^{U}_{0}\rbra{
		\int_{0}^{\infty}e^{-q\eta_{A^c}(t)}1_{\{t < \kappa_{A}\}}dt}
		n^{U}\rbra{\int_{0}^{T_0 \land \tau^{+}_{a}\land\tau^{-}_{b}}e^{-qt}f(U_{t})dt
		},
		\label{704x}
		\end{align}
	the first factor of which equals to 
		\begin{align}
		\int_{0}^{ \infty}e^{-t n^{U}\rbra{ 1 - e^{-qT_0}; A^c}}
			e^{-tn^{U}(A)}dt
		=&\frac{1}{n^{U}\rbra{ 1 - e^{-qT_0}1_{\{A^c\}}}}.
		\end{align} 
	When $X$ has bounded variation paths, we have
		\begin{align}
		&\bE^{U}_{0}\rbra{\int_{0}^{\tau^{+}_{a} \land\tau^{-}_{b}} e^{-qt}  f(U_{t} ) dt}
		\nonumber
		\\
		=&\sum_{n=0}^{\infty}
			\bE^{U}_{0}\rbra{\int_{T_0^{(n)}}^{T_0^{(n+1)} \land \tau^{+}_{a} \land\tau^{-}_{b}} 
			e^{-qt}  f(U_{t} ) dt;T_0^{(n)}< \tau^+_a \land \tau^-_b }
		\\
		=&\sum_{n=0}^{\infty}\bE^{U}_0\rbra{  e^{-qT_0}; \tau^+_a = \infty, 
				\tau^-_b = \infty}^{n}
			\bE^{U}_{0}\rbra{\int_{0}^{T_0 \land \tau^{+}_{a} \land\tau^{-}_{b}} 
			e^{-qt}  f(U_{t} ) dt}
		\label{614e}\\
		=&\frac{\bE^{U^0}_{0}\rbra{\int_{0}^{T_0 \land \tau^{+}_{a} \land\tau^{-}_{b}} 
			e^{-qt}  f(U_{t} ) dt}}{1- \bE^{U^0}_0\rbra{  e^{-qT_0}; \tau^+_a = \infty, 
				\tau^-_b = \infty}}, 
		\label{714i}
		\end{align}
	where we used the notation of the proof of Theorem \ref{ref potential density}. 
	Since $n^U=\delta_X\bE_0^U$, we obtain 
\begin{align}
	\bE^{U}_{0}\rbra{\int_{0}^{\tau^{+}_{a} \land\tau^{-}_{b}} e^{-qt}  f(U_{t} ) dt}
		=\frac{n^{U}\rbra{\int_{0}^{T_0 \land \tau^{+}_{a} \land\tau^{-}_{b}} 
			e^{-qt}  f(U_{t} ) dt}}{n^{U}\rbra{1-  e^{-qT_0}; \tau^+_a < \infty, 
				\tau^-_b = \infty}}, 
	\end{align}
	which has the same form as in the case of unbounded variation. 
	The denominator has already computed in \eqref{436}. 
	Let us compute the numerator. 
	In the case $f = 1_{(a', a]}$ for $0<a'<a$, we have
		\begin{align}
		n^{U}\rbra{\int_{0}^{ T_0 \land \tau^{+}_{a}\land\tau^{-}_{b}}
			e^{-qt}1_{\{U_{t}\in (a', a]\}}dt} %\nonumber
		=&n^{X}\rbra{ e^{-q\tau^{+}_{a'}}; \tau^+_{a'} < \infty }
			\bE^{X}_{a'}\rbra{\int_{0}^{ \tau^{+}_{a}\land\tau^{-}_{0}}
			e^{-qt}1_{\{X_{t}\in (a', a]\}}dt}
		\label{514}
		\\
		=&\frac{1}{W_{X}^{(q)}(a)}\int_{(a', a]}W_{X}^{(q)}(a-y)dy
		\label{515}
		\end{align}
	where 
	in \eqref{515} we used Theorem \ref{equation lemma} 
	and \eqref{killed potential density}. 
	Thus we obtain \eqref{701z} for $x=0$ and $y \in (0, a]$. 
	In the case $f = 1_{[b, b')}$ for $b<b' < 0$, we have
		\begin{align}
		n^{U}\rbra{\int_{0}^{ T_0 \land\tau^{+}_{a}\land\tau^{-}_{b}}
			e^{-qt}1_{\{U_{t}\in [b, b')\}}dt} 
		=&n^{X}\rbra{ e^{-q\tau^{-}_{0}}\bE^{Y}_{X(\tau^{-}_{0})}
		\rbra{ \int_{0}^{ \tau^{-}_{b}\land T_0}
			e^{-qt}1_{\{Y_{t}\in [b, b')\}}dt}
			; \tau^{-}_{0} < \tau^{+}_{a}}.
		\label{518}
		\end{align}
	Using 
	Lemma \ref{excursion Gerber-Shiu under a}, 
	we have 
		\begin{align}
		\eqref{518}
	%	=&\int \rbra{
	%		\int_{b}^{ b'} 
	%		\left(\frac{W_{Y}^{(q)}(u - b)W_{Y}^{(q)}(-y)}{W_{Y}^{(q)}(-b)}
	%		-W_{Y}^{(q)}(u-y)\right)dy}\frac{W_{X}^{(q)}(a-v)}{W_{X}^{(q)}(a)}
	%		\tilde{\Pi}_{X}(du~ dv)
	%	\label{519}
	%	\\
		=&\int_{b}^{ b'} \rbra{\int 
			\overline{\underline{r}}_Y^{(q;b,0)}(u, y)
			\frac{W_{X}^{(q)}(a-v)}{W_{X}^{(q)}(a)}
			\tilde{\Pi}_{X}(du~ dv)}dy.
		\label{520}
		\end{align}
	Using \eqref{436}, \eqref{515} and \eqref{520}, we obtain \eqref{701z} for 
	$x=0$ and $y \in [b, 0)$.\par
	{\bf{Step.2}} We calculate in the case $x < 0$.
	We have
		\begin{align}
				&\bE^{U}_{x}\rbra{\int_{0}^{ \tau^{+}_{a}\land\tau^{-}_{b}} e^{-qt}  f(U_{t})dt}
		\nonumber
		\\
		=&\bE^{Y}_{x}\rbra{
			\int_{0}^{ \tau^{+}_{0}\land\tau^{-}_{b}}e^{-qt}f(Y_t)dt}
			+\bE^{U}_{x}\rbra{ e^{-q\tau^{+}_{0}}\rbra{ 
			\int_{0}^{ \tau^{+}_{a}\land\tau^{-}_{b}}e^{-qt}f(U_t)dt}
			\circ \theta_{\tau^{+}_{0}};\tau^{+}_{0} < \tau^{-}_{b}}.
		\label{524}
		\end{align}
	Using \eqref{potential density}, 
	 we have that the first term equals to
		\begin{align}
		\int_{b}^{ 0}f(y)
		\overline{\underline{r}}_Y^{(q;b,0)}(x, y)dy.
		\label{527}
		\end{align}
	Using \eqref{exit problem} and 
	{\bf{Step.1}}, we have that the second term equals to 
		\begin{align}
		&\bE^Y_x\rbra{e^{-q\tau^+_0};\tau^+_0 < \tau^-_b}
		\bE^U_0\rbra{\int_0^{\tau^+_a \land \tau^-_b} e^{-qt}f(U_t)dt}
		\\
			=&\frac{W_U^{(q)}( x, b)}{W_U^{(q)}(  a, b)}		
		{\Bigg(}\int_0^a f(y)W_{X}^{(q)}(a-y)dy\\
			&+\int_b^0 f(y)\rbra{\int W_X^{(q)}(a-v)
			\times\overline{\underline{r}}_Y^{(q;b,0)}(u, y)
			\tilde{\Pi}_{X}(du~ dv)}dy.
			{\Bigg)}
		\label{531}
		\end{align}
	Using \eqref{527} and \eqref{531}, we obtain \eqref{701z} for $x<0$. \par
	{\bf{Step.3}} We calculate in the case $x > 0$.
	We have 
		\begin{align}
		&\bE^{U}_{x}\rbra{
			\int_{0}^{ \tau^{+}_{a}\land\tau^{-}_{b}}e^{-qt}f(U_t)dt
			}\nonumber
		\\
		=&\bE^{X}_{x}\rbra{
			\int_{0}^{ \tau^{-}_{0}\land\tau^{+}_{a}}e^{-qt}f(X_t)dt}
			+\bE^{U}_x\rbra{ e^{-q\tau^{-}_{0}}\rbra{
			\int_{0}^{ \tau^{+}_{a}\land\tau^{-}_{b}}e^{-qt}f(U_t)dt}
			\circ \theta_{\tau^{-}_{0}};\tau^{-}_{0} < \tau^{+}_{a}}.
		\label{535}
		\end{align}
	Using \eqref{killed potential density}, we have 
		the first term equals to 
		\begin{align}
		\int_{0}^{ a}f(y)
		\overline{\underline{r}}_X^{(q;0,a)}(x, y)dy.
		\label{536}
		\end{align}
	The second term equals to 
	\begin{align}
		&\bE^{X}_{x}\rbra{ e^{-q\tau^{-}_{0}}
			\bE^{U}_{X(\tau^{-}_{0})}\rbra{
			\int_{0}^{ \tau^{+}_{a}\land\tau^{-}_{b}}e^{-qt}f(U_t)dt
			};\tau^{-}_{0} < \tau^{+}_{a}}
		%\label{539}
		\nonumber
		\\
		=&
		\int
		\bE^{U}_{u}\rbra{
			\int_{0}^{ \tau^{+}_{a}\land\tau^{-}_{b}}e^{-qt}f(U_t)dt
			}
		\overline{\underline{r}}_X^{(q;0,a)}(x, v)
		\tilde{\Pi}_{X}(du~ dv),
	\label{642c}
		\end{align}
	where in \eqref{642c} we used {$(ii)$ of Theorem \ref{Gerber-Shiu}}. 
	If $f$ is $0$ on $(-\infty, 0]$, we have 
		\begin{align}
		\eqref{642c}
		=&\rbra{\frac{W_U^{(q)}(x, b)}{W_U^{(q)}(a, b)}- \frac{W_{X}^{(q)}(x)}{W_{X}^{(q)}(a)}}
			\int_0^{\infty}f(y)W_X^{(q)}(a-y)dy.
		\label{644c}
		\end{align}
	From \eqref{536} and \eqref{644c} we obtain \eqref{701z} for 
	$x>0$ and $y \in (0, a]$. 
	If $f$ is $0$ on $(0, \infty)$, we have
		\begin{align}
		\eqref{642c}
		=\int_0^\infty f(y)\rbra{\frac{W_U^{(q)}( x, b)}{W_U^{(q)}(a, b)}W_U^{(q)}(a, y)
			-W_U^{(q)}( x, y)}dy.
		\end{align}
	Thus we obtain \eqref{701z} for $x>0$ and $y\in [b, 0)$.
	}
%	Let $f$ be a non-negative measurable function.
%	Using the monotone convergence theorem twice, we have
%		\begin{align}
%		\int_\bR r_U^{(q)}(x, y)f(y)dy 
%		=&\bE^{U}_{x}\rbra{\int_0^{\infty}e^{-qt}f(U_t)dt }\\
%		=&\lim_{b \downarrow -\infty} \lim_{a\uparrow \infty}
%			\bE^{U}_{x}\rbra{\int_0^{\tau^{+}_a \land \tau^-_b}e^{-qt}f(U_t)dt }\\
%		=&\lim_{b \downarrow -\infty} 
%			\lim_{a\uparrow \infty}\int_\bR  
%			\overline{\underline{r}}_U^{(q;b,a)}( x, y)f(y)dy \\ 
%		=&\int_\bR \lim_{b \downarrow -\infty} 
%			\lim_{a\uparrow \infty} \overline{\underline{r}}_U^{(q; b, a)}( x, y)f(y)dy.
%		\label{652d}
%		\end{align}
%	Thus it follows that
%		\begin{align}
%		r_U^{(q)}(x, y)=\lim_{b \downarrow -\infty} \lim_{a\uparrow \infty} 
%		\overline{\underline{r}}_U^{(q;b,a)}(x, y)
%		\end{align}
%	for all $x, y \in \bR$.
%	We can thus obtain potential densities and one-sided killed potential densities 
%	using Theorem \ref{killed ref potential density}.

\section{Approximation problem}\label{sec app}

Let $Z$ be a spectrally negative L\'evy process. 
Let $\Psi_Z$ denote the Laplace exponent represented by \eqref{Laplace exponent}. 
For $n\in \bN$, we define 
		\begin{align}
		\Psi_{Z^{(n)}} (q) &=  \gamma_Z q  
		-\sigma_{Z}^{2} n^{2}\rbra{ 1 - e^{ q(-\frac{1}{n})} 
			+ q \rbra{ -\frac{1}{n}} } \nonumber\\
		&~~~~~~~~- \int_{\rbra{-\infty , -\frac{1}{n}}}\rbra{ 1 - e^{q y} + 
		q y 1_{\rbra{-1, -\frac{1}{n}}} (y)} \Pi_Z (dy)
		\\
		&=\delta_{Z^{(n)}}q-\int_{(-\infty, 0)}\rbra{1-e^{qy}}\Pi_{Z^{(n)}}(dy)
		\end{align}
	where 
	\begin{align}
	\delta_{Z^{(n)}}&=\gamma_Z+\sigma_{Z}^2 n + \int_{(-1, -\frac{1}{n})}(-y)\Pi_Z(dy)
	\\
	\Pi_{Z^{(n)}} &=1_{(-\infty, -\frac{1}{n})}\Pi_Z +\sigma_Z^2 n^2 \delta_{(-\frac{1}{n})}.
	\end{align}
	If we denote by $Z^{(n)}$ a L\'evy process 
	with Laplace exponent $\Psi_{Z^{(n)}}$, it is actually a compound Poisson 
	process with positive drift. 
We note that $\Psi_{Z^{(n)}}(q)\rightarrow \Psi_Z(q)$ 
for all $q \geq 0$, 
so that we have $Z^{(n)} \rightarrow Z$ in law on $\bD$. 
More precisely, by Bertoin \cite[pp.210]{Ber}, we see that there exists a coupling of 
${Z^{(n)}}$'s such that $Z^{(n)} \rightarrow Z$ uniformly on compact intervals 
almost surely, which we will call the {\em{uniformly convergent coupling}}. 
\par
	Let $X$ and $Y$ be spectrally negative L\'evy processes and suppose that 
	$X$ has unbounded variation paths and no Gaussian component. 
	For each $n\in \bN$, let $X^{(n)}$ and $Y^{(n)}$ be independent L\'evy processes 
	with Laplace exponents $\Psi_{X^{(n)}}$ and $\Psi_{Y^{(n)}}$, respectively. 
	Let $U^{(n)}$ be defined as a unique strong solution of the stochastic differential 
	equation 
	\begin{align}
	U^{(n)}_t = U^{(n)}_0 + \int_{(0, t]} 1_{\{U^{(n)}_{s-}  \geq 0\}} d X^{(n)}_{s}
	+ \int_{(0, t]} 1_{\{U^{(n)}_{s-} < 0\} } d Y^{(n)}_{s}.
	\end{align}
	\begin{Thm}
	\label{distribution convergence}
	$(U^{(n)}, \bP^{U^{(n)}}_x)$ converges in 
	distribution to $(U, \bP^U_x)$ for all $x\in\bR$. 
	\end{Thm}
	We postpone the proof of Theorem \ref{distribution convergence} until the 
	proof of Theorem \ref{resolvent convergence}.
	\begin{Rem} 
	We may expect 
	\begin{align}
	\delta_{X^{(n)}}\bP^{X^{(n)0}} \rightarrow n^X
	\text{ and } \delta_{X^{(n)}}\bP^{U^{(n)0}} \rightarrow n^U. 
	\end{align}
	The precise statements are as follows: 
	For all bounded continuous function $f$, we have
		\begin{align}
		\delta_{X^{(n)}}\bE^{X^{(n)0}}_0\rbra{\int_0^{T_0}e^{-qt}f(X^{(n)}_t)dt}
		\rightarrow
		n^X\rbra{\int_0^{T_0} e^{-qt}f(X_t)dt}
		\text{ as } n\uparrow \infty
		\label{842j}
		\end{align}
	and
		\begin{align}
		\delta_{X^{(n)}}\bE^{U^{(n)0}}_0\rbra{\int_0^{T_0}e^{-qt}f(U^{(n)}_t)dt}
		\rightarrow
		n^U\rbra{\int_0^{T_0} e^{-qt}f(U_t)dt}
		\text{ as } n\uparrow \infty.
		\label{843j}
		\end{align}
	The proofs of these formulas are straightforward, so we omit it. 
	\end{Rem}
	\begin{Lem}\label{Y convergence}
	For all non-positive $x^{(n)}$ and $x$ satisfying $x^{(n)} \rightarrow x$ 
	as $n\uparrow \infty$ and for all 
	$q>0$ and bounded continuous function $f$, we have 
	\begin{align}
	R_{Y^{(n) 0}}^{(q)} f(x^{(n)}) \rightarrow R_{Y^{0}}^{(q)} f(x) \text{ as } 
	n\uparrow \infty. 
	\label{804aa}
	\end{align}
	\end{Lem}
		\Proof{
	Using the strong Markov property, we have 
	\begin{align}
	R_{Y^{(n) 0}}^{(q)} f(x^{(n)}) 
	= R_{Y^{(n)}}^{(q)} f(x^{(n)}) 
		- \bE^{Y^{(n)}}_{x^{(n)}}(e^{-q\tau^+_0}) R_{Y^{(n) }}^{(q)} f(0)
	%= R_{Y^{(n)}}^{(q)} f(x^{(n)}) + e^{-\Phi_{Y^{(n)}}(q)x^{(n)}} R_{Y^{(n)}}^{(q)} f(0). 
	\label{805aa}
	\end{align}
	and a similar identity for $R_{Y^0}^{(q)}f(x)$. 
	Using the uniformly convergent coupling and 
	the dominated convergence theorem, we have 
	$R_{Y^{(n)}}^{(q)}f(x^{(n)}) = R_{x^{(n)}+Y^{(n)}}^{(q)}f(0)
	\rightarrow R_{x+Y}^{(q)}f(0)=R_{Y}^{(q)}f(x)$.
		Since $\Psi_{Y^{(n)}} \rightarrow \Psi_{Y}$ pointwise as $n\uparrow \infty$,
	we have $\Phi_{Y^{(n)}} \rightarrow \Phi_{Y}$ pointwise as $n\uparrow \infty$ 
	and thus 
		\begin{align}
	\bE^{Y^{(n)}}_{x^{(n)}} \rbra{ e^{-q\tau^+_0}} = e^{-\Phi_{Y^{(n)}}(q)x^{(n)}}
	\rightarrow
	e^{-\Phi_{Y}(q)x} =\bE^{Y}_x \rbra{ e^{-q\tau^+_0}}.
	%\text{ uniformly in } x \in [k, 0]. 
	\end{align}
	Thus we obtain \eqref{804aa}. 
	}
	\begin{Thm}\label{convergence}
	For all $x \in \bR$, $q>0$ and bounded continuous function $f$, we have
		\begin{align}\label{potential convergence}
		R_{U^{(n)}}^{(q)} f(x) \rightarrow R_{U}^{(q)} f (x)~as ~ n\uparrow \infty.
		\end{align}
%	For all $ q > 0$ and $f \in C_0$, we have
%		\begin{align}\label{potential convergence}
%		R_{U^{(n)}}^{(q)} f \rightarrow R_{U}^{(q)} f ~as ~ n\uparrow \infty
%		\end{align}
%	in sup norm.
%	Furthermore, we have 
%		\begin{align}\label{law convergence}
%		U^{(n)} \rightarrow U ~~~in ~law ~on~ D.
%		\end{align}
	\end{Thm}
	\Proof{
	We may assume without loss of generality that $0\leq f \leq 1$.
	We write 
		$\rho_X:= \inf_{n\in\bN}\Phi_{X^{(n)}}(q)$ and 
		$\rho_Y:= \inf_{n\in\bN}\Phi_{Y^{(n)}}(q)$. 
	Since $\Phi_Z(q)$ is strictly positive for all spectrally negative L\'evy process $Z$, 
	we have $\rho_X$ and $\rho_Y$ are strictly positive. 
	\par
	We prove \eqref{potential convergence} for $x=0$. 
	By \eqref{507v} and \eqref{506v} of Theorem \ref{ref potential density}, 
	it is sufficient to prove 
		\begin{align}
		\int_0^{\infty} e^{-\Phi_{X^{(n)}}(q)y}f(y)dy \rightarrow 
		\int_0^{\infty} e^{-\Phi_{X}(q)y}f(y)dy
		\label{705}
		\end{align}
	and
		\begin{align}
		\int R_{Y^{(n)0}}^{(q)}f(u)K_{X^{(n)}}^{(q)}(du~ dv)
		\rightarrow
		\int R_{Y^{0}}^{(q)}f(u)K_{X}^{(q)}(du~ dv) .
		\label{706}
		\end{align}
%	Since $\Psi_{X^{(n)}} \rightarrow \Psi_{X}$ as $n\uparrow \infty$,
%	we have $\Phi_{X^{(n)}} \rightarrow \Phi_{X}$ as $n\uparrow \infty$. 
	Using $\Phi_{X^{(n)}} \rightarrow \Phi_{X}$ and 
	the dominated convergence theorem, we have \eqref{705}.
	Let us prove \eqref{706}. 
	Using \eqref{314aa} with $c=1$ and changing variables, we have 
		\begin{align}
		&\int R_{Y^{(n)0}}^{(q)}f(u)K_{X^{(n)}}^{(q)}(du~ dv)
		\nonumber\\
		=&\int_{(-\infty, 0)} \Pi_{X}(du )1_{(u<-\frac{1}{n})}\int_{0}^{-u}
			e^{-\Phi_{X^{(n)}}(q)v}R_{Y^{(n)0}}^{(q)}f(u+v)
			dv
		\label{712}
		\end{align}
	and a similar identity for $(Y^0, X)$.
	We have 
		\begin{align}
		&\absol{1_{(u<-\frac{1}{n})}\int_{0}^{-u}
			e^{-\Phi_{X^{(n)}}(q)v}R_{Y^{(n)0}}^{(q)}f(u+v)
			dv} \nonumber \\
		\leq&
		\int_0^{-u}e^{-\rho_X v} \bE^{Y^{(n)}}_{u+v}\rbra{
		\int_0^{T_0} e^{-qt}dt} dv
		\\
		\leq&
		\frac{1}{q}\int_0^{-u}e^{-\rho_X v} \rbra{1-e^{\Phi_{Y^{(n)}}(q)(u+v)}}dv
		\\
		\leq&
		\frac{1}{q}\int_0^{-u}e^{-\rho_X v}\rbra{1-e^{\rho_Y(u+v)}}
		dv
		\label{715}
		\\
		\leq&\frac{1}{q\rho_X}\rbra{1-e^{\rho_X u}}\rbra{1-e^{\rho_Y u}} 
		\in L^1 \rbra{\Pi_X}. 
		\end{align}
		Thus we may apply the dominated convergence theorem to obtain 
		\begin{align}
		\lim_{n \uparrow \infty}\eqref{712}
		=\int_{(-\infty, 0)} \Pi_{X}(du )\int_{0}^{-u}
			e^{-\Phi_{X}(q)v}\rbra{R_{Y^{0}}^{(q)}f(u+v)}dv,
		\end{align}
		which shows \eqref{706}. 
		{Thus we obtain \eqref{potential convergence} for $x=0$. }\par
	For $x < 0$, \eqref{potential convergence} is obvious 
	by \eqref{504h} of Theorem \ref{ref potential density} 
		and Lemma \ref{Y convergence}. 
	\par
	We prove \eqref{potential convergence} for $x > 0$.
	By \eqref{505h} of Theorem \ref{ref potential density}, it suffices to prove 
	\begin{align}
	\underline{R}_{X^{(n)}}^{(q; 0)}f(x) \rightarrow \underline{R}_{X}^{(q; 0)}f(x)
	~~~as~n \uparrow \infty
	\label{825h}
	\end{align}
	and
	\begin{align}
	\bE^{X^{(n)}}_x\rbra{e^{-q\tau^{-}_0} R_{U^{(n)}}^{(q)}f\rbra{X^{(n)}_{\tau^-_0}} } \rightarrow
	\bE^X_x\rbra{e^{-q\tau^-_0} R_{U}^{(q)}f\rbra{X_{\tau^-_0}}}
	~as~n \uparrow \infty.
	\label{826h}
	\end{align}
	{Note that $e^{-q\tau^-_0}=e^{-q\tau^-_0}1_{(\tau^-_0 < \infty)}$ a.s. 
	Since $X$ has no Gaussian component, 
	we have
	\begin{align} 
	\inf_{t\in[0, \tau^-_0(X))}X_{t}>0 \text{ and } X_{\tau^-_0(X)}< 0 \text{ a.s. on } 
\{\tau^-_0 (X)< \infty\}. 
	\end{align}}
	For {almost every sample path with $\tau^-_0(X)<\infty$ based on} the uniformly convergent coupling 
	{of \cite[pp.210]{Ber},   
	we have 
	\begin{align}
	\inf_{t\in[0, \tau^-_0(X))}X^{(n)}_{t}\underset{n\uparrow \infty}{\rightarrow}\inf_{t\in[0, \tau^-_0(X))}X_{t} ~~
	\text{ and }~~X^{(n)}_{\tau^-_0(X)}\underset{n\uparrow \infty}{\rightarrow}X_{\tau^-_0(X)}, 
	\end{align}
	so that} 
	we have 
	\begin{align}
	\tau^-_0 {(X)}= \tau^-_0 (X^{(n)}) \text{ for large } n.
	\end{align}  
	{Therefore} we have \eqref{825h} by the dominated convergence theorem. 
	By the strong Markov property, we have
	\begin{align}
	&\bE^{X^{(n)}}_x\rbra{ e^{-q\tau^{-}_0}
			R_{U^{(n)}}^{(q)}f\rbra{X^{(n)}_{\tau^{-}_0}}}\nonumber \\
	&~~~~~~~~~~~~=\bE^{X^{(n)}}_x\rbra{
		e^{-q\tau^-_0}R_{Y^{0(n)}}^{(q)}f(X^{(n)}_{\tau^{-}_0})}
		+\bE_x^{X^{(n)}}\rbra{e^{\Phi_{Y^{(n)}}(q){X^{(n)}}(\tau^{-}_0)-q\tau^-_0}
			}R_{U^{(n)}}^{(q)}f(0). 
	\label{833h}
	\end{align}
	For the first term we have 
	\begin{align}
	\lim_{n \uparrow \infty}\bE^{X^{(n)}}_x
			\rbra{e^{-q\tau^-_0}R_{Y^{0(n)}}^{(q)}f(X^{(n)}_{\tau^{-}_0})} 
	= \bE^X_x\rbra{
		e^{-q\tau^-_0}R_{Y^{0}}^{(q)}f(X_{\tau^{-}_0})}
	\end{align}
	where we used the dominated convergence theorem 
	and Lemma \ref{Y convergence}. 
	{For the second term we have 
\begin{align}
\lim_{n\uparrow \infty}\bE_x^{X^{(n)}}\rbra{e^{\Phi_{Y^{(n)}}(q){X^{(n)}}(\tau^{-}_0)-q\tau^-_0}}R_{U^{(n)}}^{(q)}f(0)
=\bE_x^{X}\rbra{e^{\Phi_{Y}(q){X}(\tau^{-}_0)-q\tau^-_0}}R_{U}^{(q)}f(0)
\end{align}
where we used $\Phi_{Y^{(n)}}\rightarrow\Phi_Y$ and 
\eqref{potential convergence} for $x=0$. 
The proof is now complete. }
	}
%	$C_0$ denotes the set of all continuous function $f$ satisfying 
%	\begin{align}
%	\lim_{x\downarrow-\infty} f(x)= \lim_{x\uparrow \infty}f(x) = 0.
%	\end{align}
%	Then $C_0$ he \emph{Banach space} with \emph{ supremum norm}. 
%	Our objection is to prove 
%%	\begin{align}
%	R_{U^{(n)}}^{(q)}f \rightarrow R_{U^{(n)}}^{(q)}f \text{ as } n\uparrow \infty
%	\end{align}
%	for all $f \in C_0$ in supremum norm, 
%	but we haven't prove complete. 
%	Thus we prove it part of the way.
%	It is sufficient that we prove when $f$ has compact support.
For a stochastic process $Z$, $t>0$, $x\in\bR$ and 
	positive or bounded measurable function $f$, 
	we define
	\begin{align}
	P^Z_t f(x):=\bE^Z_x\rbra{f(Z_t)}.
	\end{align}
	\begin{Thm}
	\label{resolvent convergence}
	For all $q>0$, $t >0$ and $f\in C_0$, we have 
	\begin{align}
	&R^{(q)}_{U^{(n)}}f \rightarrow R^{(q)}_{U}f \text{ uniformly as }n\uparrow \infty, 
	\label{uniformly} \\ 
		&P^{U^{(n)}}_{t}f \rightarrow P^{U}_{t}f \text{ uniformly as }n\uparrow \infty. 
	\label{8119o}
	\end{align}
	\end{Thm}	
	\Proof[Proof of Theorems \ref{distribution convergence} and 
		\ref{resolvent convergence}]{
		From \eqref{uniformly} we can derive \eqref{8119o} by 
		using Theorem \ref{Feller property} 
		and \cite[Theorem $3.4.2$]{Paz}. 
		Using \cite[Theorem $19.25$]{Kal}, we can conclude that 
		$(U^{(n)}, \bP^{U^{(n)}}_x)$ converges in 
	distribution to $(U, \bP^U_x)$ for all $x\in\bR$.
	\par
	Let us prove \eqref{uniformly}. 
	We divide the proof of \eqref{uniformly} into three steps.\par
	{\bf{Step.1}} Let $k>0$ be a constant. 
	We prove $\{ \overline{W}_{U^{(n)}}^{(q)}( x)\}_{n\in\bN}$ is equicontinuous 
	in $x\in[-k, k]$. 
	For this, we prove pointwise convergence 
	$\lim_{n \uparrow \infty }\overline{W}_{U^{(n)}}^{(q)}( x)
		=\overline{W}_{U}^{(q)}( x)$. 
	Since $\{\overline{W}_{U^{(n)}}^{(q)}\}_{n \in \bN}$ is increasing and 
	continuous by Corollary \ref{refracted one-sided exit problem}, 
	the pointwise convergence implies uniform convergence in $x \in [-k, k]$, 
	thus $\{ \overline{W}_{U^{(n)}}^{(q)}( x)\}_{n\in\bN}$ is equicontinuous 
	in $x\in[-k, k]$. 
	The desired convergence is obvious for $x\leq 0$ by the definition of 
	$\overline{W}_U^{(q)}(x)$. 
	\par
	For $x>0$, it suffices to show 
	\begin{align}
	\lim_{n \uparrow \infty}\bE^{U^{(n)}}_0\rbra{e^{-q\tau^+_x}}
	=\bE^{U}_0\rbra{e^{-q\tau^+_x}}
	\label{841l}
	\end{align} 
	by Corollary \ref{refracted one-sided exit problem}. 
	By the strong Markov property, we have 
	\begin{align}
	R_{U^{(n)}}^{(q)}1_{(-\infty, x) }(0)
	=&\frac{1}{q}\rbra{1-\bE^{U^{(n)}}_0(e^{-q\tau^+_x})}
		+\bE^{U^{(n)}}_0(e^{-q\tau^+_x})
			R_{U^{(n)}}^{(q)}1_{(-\infty, x)}(x). \label{843l}
	\end{align}
	As $f^- :=1_{(-\infty, x)}$ is not continuous, we take bounded continuous 
	functions such that $f^-_m$ and $f^+_m$ such that $f^-_m \uparrow f^-$ and 
	$f^+_m \downarrow f^+ := 1_{(-\infty, x]}$. 
	Using Theorem \ref{convergence}, we have 
	$R_{U^{(n)}}^{(q)}f^\pm_m \rightarrow R_{U}^{(q)}f^\pm_m$. 
	It is obvious that $R_{U^{(n)}}^{(q)}f^\pm \rightarrow R_{U}^{(q)}f^\pm$. 
	Thus we obtain \eqref{841l}.
	\par
	{\bf{Step.2}} 
	We may assume without loss of generality that $\norm{f} =1$. 
	Let us prove 
	\begin{align}
	R^{(q)}_{U^{(n)}}f(x) \rightarrow R^{(q)}_{U}f(x)
	\text{ uniformly in }x\in[-k, k]. 
	\end{align}
	Since we have the pointwise convergence 
	by Theorem \ref{convergence}, 
	it is sufficient to prove $\{R_{U^{(n)}}^{(q)}f\}_{n\in\bN}$ is equicontinuous. 
	For all $x$, $ y\in\bR$ with $x<y$, making a computation similar to 
	{\bf{1)}} of the proof of 
	Theorem \ref{Feller property}, 
	we have 
	\begin{align}
	\absol{R_{U^{(n)}}^{(q)} f(y) - R_{U^{(n)}}^{(q)} f(x)}
	\leq 
	\frac{2}{q}\norm{f}\rbra{1- \frac{\overline{W}_{U^{(n)}}^{(q)}( x)}
				{\overline{W}_{U^{(n)}}^{(q)}( y)}}.
	\label{845k}
	\end{align}
	Let $\epsilon>0$ be a constant. 
	By {\bf{Step.1}} and since $\inf_{n\in\bN}\overline{W}_{U^{(n)}}^{(q)}(-k)=
	\inf_{n\in\bN}e^{-\Phi_{Y^{(n)}}(q)k}>0$, 
	we see that there exists $\xi >0$ such that for all $x,y \in [-k, k]$ 
	with $ 0<y-x<\xi$
	\begin{align} 
	\sup_{n\in\bN}\absol{\overline{W}_{U^{(n)}}^{(q)}( y)-\overline{W}_{U^{(n)}}^{(q)}( x)}
	\leq\epsilon\inf_{n\in\bN}\overline{W}_{U^{(n)}}^{(q)}( -k).
	\end{align}
	Then we have 
	\begin{align}
	\eqref{845k}
	\leq\frac{2}{q}\norm{f}
		\frac{
			\epsilon\inf_{n\in\bN}\overline{W}_{U^{(n)}}^{(q)}( -k)}
			{\overline{W}_{U^{(n)}}^{(q)}( y)}
	\leq\frac{2}{q}\norm{f}\epsilon , 
	\end{align}
	where we used the fact that $\overline{W}_{U^{(n)}}^{(q)}$ is increasing. 
	Therefore we conclude that $\{R_{U^{(n)}}^{(q)}f\}_{n\in\bN}$ is equicontinuous. \par
	{\bf{Step.3}} %Let $\epsilon>0$ be a constant and 
	%$\gamma=\sup_{n\in \bN}\Phi_{Y^{(n)}}(q)$.
	We prove that for any $\epsilon>0$ there is $k>0$ such that 
	\begin{align}
	\sup_{x\in (-\infty, -k)\cup(k, \infty)}
		\sup_{n\in \bN}\absol{R_{U^{(n)}}^{(q)}f(x)} < \epsilon.
	\label{855l}
	\end{align}
	For all $x<y<0$ we have 
	\begin{align}
	\absol{R_{U^{(n)}}^{(q)}f(x)}
	=&\absol{\bE_x^{U^{(n)}}\rbra{\int_0^{\tau^+_y} e^{-qt}f(U^{(n)}_t)dt}
	+\bE^{U^{(n)}}_x\rbra{e^{-q\tau^+_y}} R_{U^{(n)}}^{(q)}f(y)}\\
	\leq&\frac{1}{q}\sup_{z<y}\absol{f(z)}
	+\frac{1}{q}\sup_{m\in\bN}\bE^{Y^{(m)}}_x\rbra{e^{-q\tau^+_y}}\norm{f}.%\\
	\label{859l}
	\end{align}
	By the same argument, for all $x>y>0$, we have
	\begin{align}
	\absol{R_{U^{(n)}}^{(q)}f(x)}
	\leq& \frac{1}{q}\sup_{z>y}\absol{f(z)}
	+\sup_{m\in \bN}\bE^{X^{(m)}}_x\rbra{e^{-q\tau^-_y}}\norm{f}.
	\label{861l}
	\end{align}
	Since $f \in C_0$, there exists $k_1>0$ such that
	\begin{align} 
	\sup_{\absol{z}>k_1}\absol{f(z)} <\frac{1}{3}q\epsilon. 
	\label{862l}
	\end{align}
	Using the uniformly convergence coupling, we have for $x>y>0$
	\begin{align}
	\lim_{n\uparrow \infty}\bE^{Y^{(n)}}_{-x}\rbra{e^{-q\tau^+_{-y}}}
	=\bE^{Y}_{-x}\rbra{e^{-q\tau^+_{-y}}}~~~\text{ and }~~~
	\lim_{n\uparrow \infty}\bE^{X^{(n)}}_x\rbra{e^{-q\tau^-_y}}
	=\bE^{X}_x\rbra{e^{-q\tau^-_y}}
	\label{maru1}
	\end{align} 
	and 
	\begin{align}
	\lim_{x\uparrow \infty}\bE^{Y}_{-x}\rbra{e^{-q\tau^+_{-y}}}=0
	~~~\text{ and }~~~
	\lim_{x\uparrow \infty}\bE^{X}_x\rbra{e^{-q\tau^-_y}} =0.
	\label{maru2}
	\end{align}
	By \eqref{maru2}, there exists $k_2 > k_1$ such that 
	\begin{align}
	\bE^{Y}_{-k_2}\rbra{e^{-q\tau^+_{-k_1}}}<\frac{\epsilon}{3\norm{f}} 
	~~~\text{ and }~~~
	\bE^{X}_{k_2}\rbra{e^{-q\tau^-_{k_1}}}<\frac{\epsilon}{3\norm{f}} 
	\label{865l}
	\end{align}
	By \eqref{maru1}, 
	there exists $N\in\bN$ such that for all $n>N$ 
	\begin{align}
	\absol{\bE^{Y^{(n)}}_{-k_2}\rbra{e^{-q\tau^+_{-k_1}}}
		-\bE^{Y}_{-k_2}\rbra{e^{-q\tau^+_{-k_1}}}}<\frac{\epsilon}{3\norm{f}}
	\label{866l}
	\end{align}
	and
	\begin{align}
	\absol{\bE^{X^{(n)}}_{k_2}\rbra{e^{-q\tau^+_{k_1}}}
		-\bE^{X}_{k_2}\rbra{e^{-q\tau^+_{k_1}}}}<\frac{\epsilon}{3\norm{f}}. 
	\label{867l}
	\end{align}
	By \eqref{maru2} again, 
	there exists $k_3>k_2$ such that for all $n\leq N$
	\begin{align}
		\bE^{Y^{(n)}}_{-k_3}\rbra{e^{-q\tau^+_{-k_1}}}<\frac{\epsilon}{3\norm{f}} 
	~~~\text{ and }~~~
	\bE^{X^{(n)}}_{k_3}\rbra{e^{-q\tau^-_{k_1}}}<\frac{\epsilon}{3\norm{f}} 
	\label{868l}
	\end{align}
	Thus we obtain 
	\begin{align}
	\sup_{n\in\bN}\bE^{Y^{(n)}}_{-k_3}\rbra{e^{-q\tau^+_{-k_1}}} 
	<\frac{2\epsilon}{3\norm{f}} 
	~~~\text{ and }~~~
	\sup_{n\in\bN}\bE^{X^{(n)}}_{k_3}\rbra{e^{-q\tau^-_{k_1}}} 
	<\frac{2\epsilon}{3\norm{f}}.
	\label{870l}
	\end{align}
	By \eqref{859l}, \eqref{861l}, \eqref{862l} and \eqref{870l}, we obtain \eqref{855l}.
	\par
	The proof is complete. 
	}
	\par
\appendix
\section{Constructing generalized a refracted process by excursions} \label{construction}
In this section, we show that we can construct from $n^U$ a right-continuous 
strong Markov processes by means of the excursion theory. 
We need the following theorem which we state without proof. 
For $t\geq 0$, we denote ${\cal{D}}_t=\sigma (\omega \mapsto \omega (s):s \leq  t)$. 
\begin{Thm}[{\cite[Theorem $2$]{Sal}}]\label{Thm302g}
Let $(Z^0, \bP^{Z^0}_x)$ be a $\bR$-valued right-continuous strong Markov process stopped at $0$. 
Suppose that a $\sigma$-finite measure $n$ on $\bD$ satisfies the following conditions: 
\begin{enumerate}
\item{$n$ is concentrated on 
	$\bD^0:=\{\omega \in \bD :  \omega (0)=0, T_0(\omega) >0, 
	\omega (t)=0 \text{ for } t \geq T_0 \}$.}\label{ia}
%\item{$n^Z \rbra{\omega: \omega (0) \neq 0 } = 0$.}
\item{$n \rbra{\bD^0}=\infty$. }\label{iia}
\item{$n\rbra{ 1 - e^{-T_0} } < \infty  $. }\label{iiia}
\item{ For all $t > 0$, $A_1 \in {\cal{D}}_t $ with 
$A_1 \subset \{ T_0 > t \}$ and $A_2 \in {\cal{B}} (\bD )$, 
\begin{align}
n\rbra{ A_1 \cap \theta^{-1}_t (A_2)   }
=\int_{A_1}\bP^{Z^0}_{\omega (t)}
\rbra{Z^0\in A_2 } n\rbra{d \omega } ,\label{A01}
\end{align}
where $\theta_t$ denotes the shift operator. }\label{iva}
\item{If a measure $n^\prime$ on $\bD$ satisfies $n \geq n^\prime \geq 0$ and 
the counterpart of Condition \eqref{iva} for $n^\prime$, 
then either $n^\prime (\bD^0)=0$ or $n^\prime (\bD^0)=\infty$. }\label{va}
\end{enumerate} 
Then there is a right-continuous strong Markov process $Z$
for which $n$ is an excursion measure 
away from $0$ and $(Z^0 , \bP^{Z^0}_x)$ is the stopped process. 
\end{Thm}
To construct the strong Markov process $U$ in Section 
\ref{General}, we need to check 
that $U^0$ is a right-continuous strong Markov process and 
that $n^U$ satisfies conditions of Theorem \ref{Thm302g}. 
\begin{Lem}\label{LemA02}
The stopped process $(U^0, \bP^{U^0}_x)$ has the Markov property. 
\end{Lem}
\Proof{
It is obvious that $(U^0 , \bP^{U^0}_x) = (Y^0 , \bP^0_x)$ for $x < 0$ satisfies the 
Markov property. 
We thus need to prove that $(U^0 , \bP^{U^0}_x)$ satisfies the Markov property for $x >0$. 
Let $A_1 \in {\cal{D}}_t$ with $A_1 \subset \{ T_0 > t \}$ and $A_2\in {\cal{B}}(\bD)$. 
We write $A=A_1 %\cap \{ T_0 > t \} 
\cap \theta^{-1}_t (A_2)   $. 
By the definition of $\bP^{U^0}_x$, we have
\begin{align}
\bP^{U^0}_x\rbra{U^0\in A  }
=&\bE^X_x\rbra{\bP^{Y^0}_{y} \rbra{
w \circ Y^0 \in A}{\biggr|}_{\tiny{\begin{subarray} xy=X(\tau^{-}_{0})\\w=(X(s))_{s<\tau^{-}_{0}}\end{subarray}}}
;  \tau^-_0 \leq t}      
\\
&+\bE^X_x\rbra{\bP^{Y^0}_{y} \rbra{
w \circ Y^0 \in A}{\biggr|}_{\tiny{\begin{subarray} xy=X(\tau^{-}_{0})\\w=(X(s))_{s<\tau^{-}_{0}}\end{subarray}}}
; t > \tau^-_0}    ,  \label{A03}
\end{align}
where $w\circ w^\prime$ denotes the concatenation 
of a path $w={(w_s)}_{s<s_0}$ of finite length 
$s_0$ and a path $w^\prime = {(w^\prime_s)}_{s \geq 0}$ of infinite length: 
\begin{align}
{(w\circ w^\prime)}_s = 
\begin{cases}
w_s~~~~~~~~~&s < s_0, \\
w^\prime_{s-s_0}   & s \geq s_0.  
\end{cases}
\end{align}
By the Markov property of $Y^0$, we have 
\begin{align}
&\bE^X_x\rbra{\bP^{Y^0}_{y} \rbra{
w\circ Y^0 \in A}{\biggr|}_{\tiny{\begin{subarray} xy=X(\tau^{-}_{0})\\w=(X(s))_{s<\tau^{-}_{0}}\end{subarray}}}
;  \tau^-_0 \leq t}  
\\
=&\bE^X_x\rbra{\bP^{Y^0}_{y} \rbra{
w \circ Y^0 \in A_1, \,
\rbra{Y^0_s}_{s \geq t-u} \in A_2}{\Biggr|}_{\tiny{\begin{subarray} xy=X_u\\w=(X(s))_{s<u}\\u=\tau^-_0\end{subarray}}}
; \tau^-_0 \leq t} \\
=&\bE^X_x\rbra{\bE^{Y^0}_y \rbra{
1_{\{w \circ Y^0 \in A_1\}}
\bP^{Y^0}_{y^\prime}\rbra{Y^0 \in A_2}\Big{|}_{y^\prime=Y^0_{t-u}}}{\Biggr|}_{\tiny{\begin{subarray} xy=X_u\\w=(X(s))_{s<u}\\u=\tau^-_0\end{subarray}}}
; \tau^-_0 \leq t}\\
=&\bE^U_x\rbra{
1_{\{U \in A_1\}}
\bP^{U^0}_{y}\rbra{U^0 \in A_2}\Big{|}_{y=U_t}
; \tau^-_0 \leq t} .
\end{align} 
We can do a similar argument for \eqref{A03}.
So we obtain
\begin{align}
\bP^{U^0}_x\rbra{U^0 \in  A  }
=\int_{A_1%\cap \{ T_0 > t \} 
}\bP^{U^0}_{\omega(t)}
\rbra{U^0\in A_2 } \bP^{U^0}_x
\rbra{U^0\in d \omega }. 
\end{align}
The proof is complete. 
}
\begin{Lem}
The stopped process $U^0$ has the strong Markov property. 
\end{Lem}
\Proof{
Fix $t >0$. 
By the proof of \cite[Theorem $1$ of Section $2.3$]{ChuWal}, 
it is sufficient to prove that 
$x \mapsto \bE^{U^0}_x\rbra{ f(U^0_t)}$ is continuous for all 
bounded continuous function $f$ with $f(0)=0$. 
Continuity at $x < 0$ is obvious, by the Feller property of $Y^0$. 
Left-continuity at $x=0$ is also obvious. 
Right-continuity at $x=0$ follows from the fact that 
$\bP^{U^0}_y\rbra{T_0 \in \cdot}\underset{y\rightarrow 0}{\rightarrow} \delta_0$.
Let us consider continuity at $x>0$. 
\begin{align}
\bE^{U^0}_y\rbra{ f(U^0_t)}   
=&\bE^{U^0}_y \rbra{f(U^0_t); \tau^-_0 \land t < T_x}
+\bE^{X}_y\rbra{ f(X_t) ; T_x \leq t <\tau^-_0 }\\
&+\bE^{X}_y\rbra{ \bE^{Y^0}_{y^\prime}\rbra{f(Y^0_{t-u})}
	{\biggr|}_{\tiny{\begin{subarray} xy^\prime=X(u)\\
u=\tau^-_0\end{subarray}}} ;T_x \leq  \tau^-_0 \leq t} . \label{A013}
\end{align}
Note that we have $\bP^X_0 \rbra{\lim_{y \rightarrow 0}T_y = 0 }=1$ 
by the assumption that $X$ is spectrally negative and of bounded variation. 
Since $X$ and $Y^0$ have c\`adl\`ag paths, 
we have the following identities: 
\begin{align}
\bE^{U^0}_y \rbra{f(U^0_t); \tau^-_0 \land t < T_x}
&\leq \norm{f} \bP^X_0 \rbra{\tau^-_{-\frac{x}{2}}<T_{x-y}}
\underset{y\rightarrow x}{\rightarrow}0, 
\\
\bE^{X}_y\rbra{ f(X_t) ; T_x \leq t <\tau^-_0 } 
&=\bE^{X}_y\rbra{ 
\bE^X_x \rbra{f(X_{t-u}); t <\tau^-_0}\Big{|}_{u=T_x}  ; T_x \leq t \land\tau^-_0 } \\
&\underset{y\rightarrow x}{\rightarrow}
\bE^X_x\rbra{  f(X_t) ; t<\tau^-_0}, 
\end{align}
\begin{align}
&\bE^{X}_y
\rbra{ \bE^{Y^0}_{y^\prime}\rbra{f(Y^0_{t-u})}
	{\biggr|}_{\tiny{\begin{subarray} xy^\prime=X(u)\\
u=\tau^-_0\end{subarray}}} ;T_x \leq  \tau^-_0 \leq t} 
\\
&=\bE^{X}_y\rbra{
\bE^X_x
\rbra{ \bE^{Y^0}_{y^\prime}\rbra{f(Y^0_{t-u-v})}
	{\biggr|}_{\tiny{\begin{subarray} xy^\prime=X(u)\\
u=\tau^-_0\end{subarray}}} ;  \tau^-_0 \leq t} 
\biggr{|}_{v=T_x}
; T_x \leq  \tau^-_0 \land t}\\
&\underset{y\rightarrow x}{\rightarrow}
\bE^{X}_x
\rbra{ \bE^{Y^0}_{y^\prime}\rbra{f(Y^0_{t-u})}
	{\biggr|}_{\tiny{\begin{subarray} xy^\prime=X(u)\\
u=\tau^-_0\end{subarray}}} ; \tau^-_0 \leq t}. 
\end{align}
The proof is now complete. 
}
\begin{Lem}\label{Lem304h}
The measure $n=n^U$ satisfies Conditions 
\eqref{ia}, \eqref{iia}, \eqref{iiia}, \eqref{iva} and 
\eqref{va} in Theorem \ref{Thm302g}. 
\end{Lem}
\Proof{ 
It is obvious by definition that $n^U$ satisfies \eqref{ia} and \eqref{iia}. 
\par
Let us prove \eqref{iiia}. 
By the definition of $n^U$ and by \eqref{one-sided exit problem}, we have 
\begin{align}
n^U\rbra{1 - e^{-T_0}} 
&= n^X\rbra{1 - e^{-\tau^-_0} \bE^{Y}_{X_{\tau^-_0}}
\rbra{e^{-\tau^+_0}}1_{\{\tau^-_0 < \infty\}}  }\\
&= n^X\rbra{1 - e^{-\tau^-_0} e^{\Phi_Y(1)X_{\tau^-_0}}1_{\{\tau^-_0 < \infty\}}  }. 
\label{A24}
\end{align}
We let $q^\prime = 1 \lor \inf \{q > 0: \Phi_X(q) > \Phi_Y (1) \}$. 
Since $n^X\rbra{ 1 - e^{-q^\prime T_0}}$ is finite, we obtain 
\begin{align}
\eqref{A24}&\leq 
n^X\rbra{1 - e^{-q^\prime\tau^-_0} 
e^{\Phi_X(q^\prime)X_{\tau^-_0}}1_{\{\tau^-_0 < \infty\}}  }\\
&=n^X\rbra{1 - e^{-q^\prime\tau^-_0} 
\bE^X_{X_{\tau^-_0}}\rbra{e^{-q^\prime\tau^+_0}}1_{\{\tau^-_0 < \infty\}}  }
<\infty.
\end{align}
\par
The proof of \eqref{iva} 
is the same as that of the Markov property of $(U^0 , \bP^{U^0}_x)$ 
for $x > 0$ in Lemma \ref{LemA02}. 
\par
Let us prove \eqref{va}. 
We define the $\sigma$-finite measure $n^{\prime\prime}$ by 
	\begin{align}
&n^{\prime\prime}\rbra{ F\rbra{(U_t)_{t < \tau^{-}_{0}}, (U_{t+\tau^{-}_{0}})_{t\geq 0}}}
=n^{\prime}\rbra{ \mathbb{E}^{X^0}_{y}
\rbra{F\rbra{w, (X^0_t)_{t\geq 0}}}
{\biggr|}_{\tiny{\begin{subarray}
xy=U(\tau^{-}_{0})\\w=(U(t))_{t<\tau^{-}_{0}}\end{subarray}}}}
	\end{align}
for all non-negative measurable functional $F$. 
Then $n^{\prime\prime}$ satisfies the Markov property for 
${\{\bP^{X^0}_x\}}_{x\in\bR\backslash \{0\}}$. 
By the definition of $n^U$, we have $n^X\geq n^{\prime\prime} \geq 0$. 
By \cite[Proposition 1]{Sal}, $n^X$ satisfies Condition \eqref{va} and 
we obtain either $n^{\prime\prime}\rbra{\bD_0}=0$ 
or $n^{\prime\prime}\rbra{\bD_0}=\infty$, which yields we have 
either $n^{\prime}\rbra{\bD_0}=0$ 
or $n^{\prime}\rbra{\bD_0}=\infty$. 
\par
The proof is complete. 
}

\noindent 
\textbf{Acknowledgments.}
The authors were supported by JSPS-MAEDI Sakura program. 
The second author was supported by MEXT KAKENHI grant no.'s 26800058
and 15H03624.

%%%%% references %%%%%
% 文献情報はMathSciNetから取得するのが便利
% http://www.ams.org/mathscinet/

\end{document}